\documentclass[a4pitaper,11pt]{amsart}
\usepackage{amsmath,amsthm,amssymb}
\usepackage[mathscr]{eucal}
 \usepackage{cite}
\usepackage{upgreek}
\usepackage[bookmarks,bookmarksnumbered]{hyperref}
\usepackage{enumerate}

\numberwithin{equation}{section}

\newcommand{\car}{\curvearrowright}

\theoremstyle{plain}
\newtheorem{main}{Theorem}
\newtheorem{mcor}[main]{Corollary}

\newtheorem{theorem}{Theorem}[section]

\newtheorem{lemma}[theorem]{Lemma}
\newtheorem{proposition}[theorem]{Proposition}
\newtheorem{corollary}[theorem]{Corollary}
\theoremstyle{definition}

\newtheorem{definition}[theorem]{Definition}
\newtheorem{example}[theorem]{Example}
\newtheorem{notation}[theorem]{Notation}
\newtheorem{remark}[theorem]{Remark}

\newtheorem{assumption}[theorem]{Assumption}
\begin{document}

\title[W$^*$-superrigidity for coinduced actions]
{W$^*$-superrigidity for coinduced actions}
\author[Daniel Drimbe]{Daniel Drimbe}

\address{Mathematics Department; University of California, San Diego, CA 90095-1555 (United States).}
\email{ddrimbe@ucsd.edu}
\thanks{The author was partially supported by NSF Career Grant DMS \#1253402.}

\begin{abstract} 

We prove W$^*$-superrigidity for a large class of coinduced actions. We prove that if $\Sigma$ is an amenable almost-malnormal subgroup of an infinite conjugagy class (icc) property (T) countable group $\Gamma$, the coinduced action $\Gamma\curvearrowright X$ from an arbitrary probability measure preserving action $\Sigma\curvearrowright X_0$ is W$^*$-superrigid. We also prove a similar statement if $\Gamma$ is an icc non-amenable group which is measure equivalent to a product of two infinite groups. In particular, we obtain that any Bernoulli action of such a group $\Gamma$ is W$^*$-superrigid.

\end{abstract}

\maketitle
\section{Introduction and statement of the main results}
\subsection{Introduction.}
To every measure preserving action $\Gamma\curvearrowright (X,\mu)$ of a countable group $\Gamma$ on a standard probability space $(X,\mu)$, one associates {\it the group measure space von Neumann algebra} $L^{\infty}(X)\rtimes\Gamma$ \cite{MvN36}. If the action $\Gamma\curvearrowright X$ is free, ergodic and probability measure preserving (pmp), then $L^{\infty}(X)\rtimes\Gamma$ is a II$_1$ factor which contains $L^\infty(X)$ as a {\it Cartan subalgebra}, i.e. a maximal abelian von Neumann algebra whose normalizer generates $L^{\infty}(X)\rtimes\Gamma$. 
The classification of group measure space II$_1$ factors $L^\infty(X)\rtimes\Gamma$ is a central problem in the theory of von Neumann algebras. Two free ergodic pmp actions $\Gamma\curvearrowright(X,\mu)$ and $\Lambda\curvearrowright (Y,\nu)$ on standard probability spaces $(X,\mu)$ and $(Y,\nu)$ are said to be W$^*${\it -equivalent} if $L^{\infty}(X)\rtimes\Gamma$ is isomorphic to $L^{\infty}(Y)\rtimes \Lambda.$

If the groups are amenable, the classification up to W$^*$-equivalency has been completed in the 1970s. More precisely, the celebrated theorem of Connes \cite{Co76} asserts that all II$_1$ factors arising from free ergodic pmp actions of countable amenable groups are isomorphic to the hyperfinite II$_1$ factor. In contrast, the non-amenable case is much more challenging and it has led to a beautiful {\it rigidity theory} in the sense that one can deduce conjugacy from W$^*$-equivalence. A major breakthrough in the classification of II$_1$  factors was made by Popa between 2001-2004 through the invention of deformation/rigidity theory (see \cite{Po07, Va10a, Io12a} for surveys). In particular, he obtained the following W$^*$-rigidity result: let $\Gamma\curvearrowright X$ be a free ergodic pmp action of an infinite conjugacy class (icc) countable group $\Gamma$ which has an infinite normal subgroup with the relative property (T) and let $\Lambda\curvearrowright Y:=Y_0^\Lambda$ be a Bernoulli action of a countable group $\Lambda.$ Popa proved that if the two actions are W$^*$-equivalent, then the actions are {\it conjugate} \cite{Po03,Po04}, i.e. there exist a group isomoprhism  $d:\Gamma\to\Lambda$ and a measure space isomorphism $\theta :X \to Y$ such that $\theta(gx)=d(g)\theta(x)$ for all $g\in\Gamma$ and almost everywhere (a.e.) $x\in X.$ 

The most extreme form of rigidity for an action $\Gamma\curvearrowright (X,\mu)$ is {\it W$^*$-superrigidity}, i.e. whenever $\Lambda\curvearrowright(Y,\nu)$ is a free ergodic pmp action W$^*$-equivalent to $\Gamma\curvearrowright(X,\mu)$, then the two actions are conjugate. 
A few years ago, Peterson was able to show the existence of {virtually} W$^*$-superrigid actions \cite{Pe09}. Soon after, Popa and Vaes discovered the first concrete families of W$^*$-superrigid actions \cite{PV09}. Ioana then proved in \cite{Io10} a general W$^*$-superrigidity result for Bernoulli actions.

{\bf Theorem} (Ioana, \cite{Io10}). If $\Gamma$ is an icc property (T) group and $(X_0,\mu_0)$ is a non-trivial standard probability space, then the Bernoulli action $\Gamma\curvearrowright (X_0,\mu_0)^\Gamma$  is W$^*$-superrigid.

The main ingredient of his proof was the discovery of a beautiful dichotomy result for abelian subalgebras of II$_1$ factors coming from Bernoulli actions.

Using a similar method, Ioana, Popa and Vaes were able to prove later that any Bernoulli action of an icc non-amenable group which is a product of two infinite groups is also W$^*$-superrigid \cite{IPV10}. A few years ago Boutonnet extended these results to Gaussian actions in \cite{Bo12}. Several other classes of W$^*$-superrigid actions have been found in \cite{FV10,CP10,HPV10,Va10b,CS11,CSU11,PV11,PV12,CIK13,CK15,Dr15,GITD16}.




\subsection{Statement of the main results.}
Our first theorem is a generalization of Ioana's W$^*$-superrigidity result \cite[Theorem A]{Io10} to coinduced actions. Before stating the theorem, we explain first the terminology that we use starting with the notion of coinduced actions (see e.g. \cite{Io06b}). 

\begin{definition}\label{coinduce}
Let $\Gamma$ be a countable group and let $\Sigma$ be a subgroup. Let $\phi:\Gamma/\Sigma \to \Gamma$ be a section. Define the cocycle $c:\Gamma\times\Gamma/\Sigma\to \Sigma$ by the formula $$c(g,i)=\phi^{-1}(gi)g\phi(i),$$ for all $g\in\Gamma$ and $i\in \Gamma/\Sigma.$\\
Let $\Sigma\overset{\sigma_{0}}{\curvearrowright} (X_0,\mu_0)$ be a pmp action, where $(X_0,\mu_0)$ is a non-trivial standard probability space. We define an action $\Gamma\overset{\sigma}{\curvearrowright} X_0^{\Gamma/\Sigma}$, called the coinduced action of $\sigma_0$, as follows: 
$$\sigma_{g}((x_{i})_{i\in \Gamma/\Sigma})=(x'_{i})_{i\in \Gamma/\Sigma}, \text{  where } x'_{i}=c(g^{-1},i)^{-1}x_{g^{-1}i}.$$
\end{definition}

Note the following remarks:
\begin{itemize}
\item $\sigma$ is a pmp action of $\Gamma$ on the standard probability space $X_0^{\Gamma/\Sigma}$.
\item if we consider the trivial action of $\Lambda=\{e\}$ on $X_0$, then the coinduced action of $\Gamma$ on $X_0^{\Gamma/\{e\}}=X_{0}^{\Gamma}$ is the Bernoulli action.
\end{itemize}

Recall that an inclusion $\Gamma_0\subset\Gamma$ of countable groups has the {\it relative property} (T) if for every $\epsilon>0$, there exist $\delta>0$ and a finite subset $F\subset\Gamma$ such that if $\pi:\Gamma\to \mathcal U (K)$ is a unitary representation and $\xi\in K$ is a unit vector satisfying $\|\pi(g)\xi-\xi\|<\delta$, for all $g\in F$, then there exists $\xi_0\in K$ such that $\|\xi-\xi_0\|<\epsilon$ and $\pi(h)\xi_0=\xi_0$, for all $h\in \Gamma_0.$  The group $\Gamma$ has the {\it property} (T) if the inclusion $\Gamma\subset \Gamma$ has the relative property (T). To give some examples, note that $\mathbb Z^2\subset \mathbb Z^2\rtimes SL_2(\mathbb Z)$ has the relative property (T) and $SL_n(\mathbb Z),$ $n\ge 3$, has the property (T) \cite{Ka67,Ma82}.

Finally, we say that a subgroup $\Sigma$ of a countable group $\Gamma$ is called {\it n-almost malnormal} if for any $g_1,g_2,...,g_n\in \Gamma$ such that $g_i^{-1}g_j\notin \Sigma$ for all $i\neq j$, the group $\cap_{i=1}^n g_i\Sigma g_i^{-1}$ if finite. The subgroup $\Sigma$ is called {\it almost malnormal} if it is n-almost malnormal for some $n\ge 1.$ 

\begin{main}\label{1}
Let $\Gamma$ be an icc group which admits an infinite normal subgroup $\Gamma_0$ with relative property (T) and let $\Sigma$ be an amenable almost malnormal subgroup of $\Gamma$. Let $\sigma_0$ be a pmp action of $\Sigma$ on a non-trivial standard probability space $(X_0,\mu_0)$ and denote by $\sigma$ the coinduced action of $\Gamma$ on $X:=X_0^{\Gamma/\Sigma}$.
Then $\Gamma\overset{\sigma}\curvearrowright X$ is $W^*$-superrigid.
\end{main}


\begin{example}
 In particular, Theorem \ref{1} can be applied for
$\Gamma=SL_3(\mathbb Z)$ and $\Sigma=\langle A\rangle$, where 
$A=
\begin{bmatrix}
0&1&1\\
-1&0&0\\
0&-1&0\\
\end{bmatrix}
$
\cite[Section 7]{PV06}. See \cite{PV06} for more concrete examples of amenable almost malnormal subgroups of $PSL_n(\mathbb Z)$, $n\ge 
3$. 
See also \cite[Theorem 1.1]{RS10}, a result which proves the existence of amenable almost malnormal subgroups of torsion-free uniform lattices in connected semisimple real algebraic groups with no compact factors.
\end{example}

We now generalize Ioana-Popa-Vaes' result \cite[Theorem 10.1]{IPV10} to coinduced actions. First, recall that two countable groups $\Gamma$ and $\Lambda$ are called {\it measure equivalent} in sense of Gromov if there exist two commuting free measure preserving actions of $\Gamma$ and $\Lambda$ on a standard measure space $(\Omega, m)$, such that the actions of $\Gamma$ and $\Lambda$ each admit a finite measure fundamental domain \cite{Gr91}. Natural examples of measure equivalent groups are provided by pairs of lattices $\Gamma,\Lambda$ in an unimodular locally compact second countable
group. 

\begin{main}\label{2}
Let $\Gamma$ be an icc non-amenable group which is measure equivalent to a product of two infinite groups. Let $\Sigma$ be an amenable almost malnormal subgroup and let $\sigma_0$ be a pmp action of $\Sigma$ on a non-trivial standard probability space $(X_0,\mu_0)$ and denote by $\sigma$ the coinduced action of $\Gamma$ on $X:=X_0^{\Gamma/\Sigma}$.\\ 
Then  $\Gamma\overset{\sigma}\curvearrowright X$ is $W^*$-superrigid.
\end{main}
 
See Theorem \ref{2'} for a more general statement in which it is assumed instead that $\Gamma$ is measure equivalent to a group $\Lambda_0$ whose group von Neumann algebra $L(\Lambda_0)$ is not prime. Note that Theorems \ref{1} and \ref{2} provide a complementary class of W$^*$-superrigid coinduced actions from the one found in \cite[Corollary 1.4]{Dr15}.

\begin{example}
 A more general statement of Theorem \ref{2} can be appplied for $\Sigma\subset\Gamma=\Delta\wr\Sigma$ with $\Delta$ non-amenable and $\Sigma$ amenable (see Remark \ref{remark.}).
\end{example}

The following remark shows that if $\Sigma$ is not almost malnormal, the action $\Gamma\car X$ is not necessary W$^*$-superrigid. To put this in context, we recall first the notion of OE-superrigidity and Singer's result \cite{Si55}. Two actions $\Gamma\car X$ and $\Lambda\car Y$ are {\it orbit equivalent} (OE) if there exists a measure space isomorphism $\theta:X\to Y$ such that $\theta(\Gamma x)=\Lambda \theta(x),$ for a.e. $x\in X.$ A pmp action $\Gamma\car X$ is {\it OE-superrigid} if whenever $\Lambda\car Y$ is a free ergodic pmp action which is OE to $\Gamma\car X$, then the two actions are conjugate.

Singer proved in \cite{Si55} that two free ergodic pmp actions $\Gamma\car X$ and $\Lambda\car Y$ are OE if and only if there exists an isomoprhism of the group measure space algebras $L^\infty(X)\rtimes\Gamma$ and $L^\infty(Y)\rtimes\Lambda$ which preserves the Cartan algebras $L^\infty(X)$ and $L^\infty(Y)$. In particular, W$^*$-superrigidity implies OE-superrigidity.

\begin{remark}
If $\Sigma$ is not almost malnormal, the action $\Gamma\car X$ may fail to be W$^*$-superrigid. Indeed, suppose $\Gamma$ is an icc group which splits as a direct product $\Gamma=\Sigma\times\Delta$, with $\Sigma$ amenable and $\Delta$ a non-amenable group. Connes and Jones have found in \cite{CJ82} a class of groups $\Sigma$ and a class of free ergodic pmp actions $\Sigma\overset{\sigma_0}{\car} X_0$ for which the coinduced action $\Gamma\car X$ of $\sigma_0$ is not W$^*$-superrigid. Precisely, they have proven that $M:=L^\infty(X)\rtimes\Gamma$  is {\it McDuff}, i.e. $M\simeq M\bar\otimes R$, where $R$ is the hyperfinite II$_1$ factor. However, \cite[Corollary 1.3]{Dr15} implies that $\Gamma\car X$ is OE-superrigid. 

\end{remark}



Note that Theorem \ref{2} extends the class of groups whose Bernoulli actions are W$^*$-superrigid. Therefore we record the following result.

\begin{mcor}\label{3}
Let $\Gamma$ be an icc non-amenable group which is measure equivalent to a product of two infinite groups. Let  $(X_0,\mu_0)$ be a non-trivial standard probability space. Then the Bernoulli action $\Gamma\curvearrowright X_0^{\Gamma}$ is $W^*$-superrigid.
\end{mcor}

We recall the well known theorem due to Borel which asserts that every connected non-compact semisimple Lie group contains a lattice (see \cite{Bo63}  and \cite[Theorem 14.1]{Ra72}). Using this, we obtain an immediate consequence of Corollary \ref{3}.

\begin{mcor}\label{4}
Let $\Gamma$ be an icc lattice in a product $G=G_1\times\dots \times G_n$ of $n\ge 2$ connected non-compact semisimple Lie groups and let $(X_0,\mu_0)$ be a non-trivial standard probability space. Then the Bernoulli action $\Gamma\curvearrowright X_0^{\Gamma}$ is $W^*$-superrigid.
\end{mcor}

Note that a combination of Popa's cocycle superrigidity theorem for product groups \cite{Po06} and the results on uniqueness of Cartan subalgebras from \cite{PV12} already proves Corollary \ref{4}, but only in the case when each factor $G_1,\dots,G_n$ is of rank one.

\subsection{Comments on the proof of Theorem \ref{2}}

For obtaining the proofs of Theorem \ref{1} and Theorem \ref{2}, we adapt the proofs used by Ioana \cite{Io10} and Ioana-Popa-Vaes \cite{IPV10} to the context of coinduced actions. We outline briefly and informally the proof of Theorem \ref{2} since it has as a consequence Corollary \ref{3}.

To this end, let $\Gamma$ be an icc group and let $\Sigma$ be an almost malnormal subgroup. Assume $\Gamma$ is measure equivalent to a product $\Lambda_0=\Lambda_1\times\Lambda_2$ of two countable groups. By \cite{Fu99}, $\Gamma$ and $\Lambda_0$ must have stably orbit equivalent actions. To simplify notation, assume there exist free ergodic pmp actions of $\Gamma$ and $\Lambda_0$ on a probability space $(Y_0,\mu)$ whose orbits are equal, almost everywhere. Thus, $L^\infty(Y_0)\rtimes\Gamma=L^\infty(Y_0)\rtimes\Lambda_0.$\\
Suppose $\Sigma\car X_0$ is a pmp action on a non-trivial standard probability space and let $\Gamma\overset{\sigma}\car X:=X_0^{\Gamma/\Sigma}$ be the corresponding coinduced action. Our goal is to show that $\Gamma\overset{\sigma}{\car} X$ is W$^*$-superrigid.
Assume that $\Lambda\car Y$ is an arbitrary free ergodic pmp action such that
$$
M:=L^\infty(X)\rtimes\Gamma=L^\infty(Y)\rtimes\Lambda.
$$

First, we reduce the problem to showing that the Cartan subalgebras $L^\infty(X)$ and $L^\infty(Y)$ are unitarily conjugated. 
We do this by proving in Section \ref{sectionr} a cocycle superrigidity theorem for $\Gamma\car X$. Combined with \cite[Theorem 5.6]{Po05}, we obtain that $\Gamma\car X$ is OE-superrigid. Therefore, by a result of Singer \cite{Si55} it is enough to show that $L^\infty(Y)$ is unitarily cojugate to $L^\infty(X)$ in $M$. We note that this is actually equivalent to $L^\infty(Y)\prec_M L^\infty(X)$, by \cite[Theorem A.1]{Po06b}. See Section \ref{intertwining} for the definition of Popa's intertwining symbol "$\prec$". 

As is \cite{Io10}, we make use of the decomposition $M=L^\infty(Y)\rtimes\Lambda$ via the comultiplication $\Delta:M\to M\bar\otimes M$ defined by $\Delta(bv_\lambda)=bv_\lambda\otimes\lambda$ , for all $b\in L^\infty(Y)$ and $\lambda\in\Lambda$, introduced in \cite{PV09}. Here we denote by $\{v_\lambda\}_{\lambda\in\Lambda}$ the canonical unitaries implementing the action of $\Lambda$ on $L^\infty(Y)$. The next step is to prove that there exists a unitary $u\in M\bar\otimes M$ such that
\begin{equation}\label{aa}
 u\Delta(L(\Gamma))u^*\subset L(\Gamma)\bar\otimes L(\Gamma).
\end{equation}
This is obtained in two steps. A main technical contribution of our paper is to use the rigidity of $\Gamma$ inherited from the product structure of $\Lambda_0$ through measure equivalence. We do this in Section \ref{sectionr} by introducing an "amplified" version of the comultiplication map $\Delta$ which is defined on the larger von Neumann algebra $(L^\infty(Y_0)\bar\otimes L^\infty(X))\rtimes\Gamma$. Combined with the spectral gap rigidity theorem for coinduced actions (Theorem \ref{converge}) proved in Section \ref{sectioni}, we obtain the conclusion \eqref{aa}.

In Section \ref{sectionio}, following Ioana's idea \cite{Io10}, we obtain a dichotomy theorem for certain abelian algebras. The result is a straightforward adaptation of \cite[Theorem 5.1]{IPV10} to coinduced actions and has two consequences. First, we obtain
$$
\Delta(L^\infty(X))'\cap (M\bar\otimes M)\prec L^\infty(X)\bar\otimes L^\infty(X).
$$
Second, it implies a weaker version of Popa's conjugacy criterion adapted to coinduced actions. This will altogether prove Theorem \ref{2}.

\subsection{Acknowledgements} 

I am very grateful to my advisor Adrian Ioana for all the help given through valuable discussions. I would also like to thank R\'emi Boutonnet for helpful comments about the paper.

\section{Preliminaries}\label{sectionp}

\subsection{Terminology.}

A von Neumann algebra $M$ is called {\it tracial} if it is equipped with a faithful normal tracial state $\tau.$ We denote by $L^2(M)$ the completion of $M$ with respect to the norm $\|x\|_2=\sqrt{\tau(x^*x)}$. For $Q\subset M$, a unital von Neumann subalgebra of $M$, we denote by $e_Q:L^2(M)\to L^2(Q)$ the orthogonal projection onto $L^2(Q).$ We denote by $E_Q:M\to Q$, the conditional expectation onto $Q$. {\it Jones' basic construction} of the inclusion $Q\subset M$ is defined as the von Neumann subalgebra of $\mathbb B (L^2(M))$ generated by $M$ and $e_Q$.

Denote by $\mathcal U(M)$ the group of unitary elements of $M$ and by $\mathcal{N}_M(Q)=\{u\in\mathcal U(M)|uQu^*=Q\}$ the {\it normalizer} of $Q$ inside $M$. Denote also by $Q'\cap M=\{x\in M|xq=qx, \text{ for all } q\in Q\}$ the relative commutant of $Q$ in $M$ and by $\mathcal Z(M)=M\cap M'$ the center of $M$.


\subsection{Popa's intertwining by bimodules.}\label{intertwining}

We recall from \cite[Theorem 2.1 and Corollary 2.3]{Po03} Popa's intertwining by bimodules tehnique which is fundamental to deformation/rigidity theory.

\begin{theorem}\label{corner}\cite{Po03}
Let $P$ and $Q$ be (not necessarily unital) subalgebras of a tracial von Neumann algebra $M$. The following are equivalent:
\begin{itemize}
\item There exist non-zero projections $p\in P, q\in Q$, a $*$-homomorphism $\varphi :pPp\to qQq$ and a non-zero partial isometry $v\in pMq$ such that $xv=v\varphi(x),$ for all $x\in pPp.$
\item For any group $\mathcal U\subset\mathcal U(P)$ such that $\mathcal U''=P$ there is no sequence $\{u_n\}\subset\mathcal U$ such that $\|E_Q(xu_ny)\|_2\to 0, $ for all $x,y\in M.$ 
\end{itemize}
If one of these conditions holds true, then we write $P\prec_M Q$ and we say that a corner of $P$ embeds into $Q$. If $Pp'\prec_M Q$ for any non-zero projection $p'\in P'\cap pMp,$ then we write $P\prec^s_M Q.$
\end{theorem}


\subsection{Bimodules and weak containment.}
Let $M,N$ be tracial von Neumann algebras. An {\it $M$-$N$-bimodule} $_M\mathcal H_N$ is a Hilbert space $\mathcal H$ equipped with two commuting normal unital $*$-homomorphisms $M\to B(\mathcal H)$ and $N^{\text{op}}\to B(\mathcal H)$. An $M-N$-bimodule $_M\mathcal H_N$ is {\it weakly contained} in a $M$-$N$-bimodule $_M\mathcal K_N$ and we write $_M\mathcal H_N \underset{weak}{\subset}$ $_M\mathcal K_N$ if for any $\epsilon>0$, finite subsets $F\subset M, G\subset N$ and $\xi\in\mathcal H$, there exist $\eta_1,\dots,\eta_n\in\mathcal K$ such that 
$$
|\langle x\xi y,\xi\rangle-\sum_{i=1}^n \langle x\eta_iy,\eta_i\rangle|\leq\epsilon,\, \text{  for all  }\, x\in F, y\in G. 
$$
Given two bimodules $_M\mathcal H_N$ and $_N\mathcal K_P$, one can define the {\it Connes tensor product} $\mathcal H\otimes_{N}\mathcal K$ which is an $M$-$P$ bimodule (see \cite[V.Appendix B]{Co94}). If $_M\mathcal H_N \underset{weak}{\subset}$ $_M\mathcal K_N$, then $_M\mathcal H\otimes_N\mathcal L_P \underset{weak}{\subset}$ $_M\mathcal K\otimes_N\mathcal L_P$, for any $N$-$P$ bimodule $\mathcal L.$

\subsection{Relative amenability.}
Let $(M,\tau)$ be a tracial von Neumann algebra. Let $p\in M$ be a projection and $P\subset pMp,Q\subset M$ be von Neumann subalgebras. Following \cite[Definition 2.2]{OP07}, we say that $P\subset pMp$ is {\it amenable relative to $Q$ inside $M$} if there exists a positive linear functional $\varphi:p\langle M,e_Q\rangle\to\mathbb C$ such that $\varphi_{|pMp}=\tau$ and $\varphi$ is $P$-central. We say that $M$ is {\it amenable} if $M$ is amenable relative to $\mathbb C1$ inside $M$.  \\
By \cite[Section 2.2]{OP07}, $P$ is amenable relative to $Q$ inside $M$ if and only if $_ML^2(Mp)_P$ is weakly contained in $_ML^2(\langle M,e_Q\rangle p)_P$.

A von Neumann subalgebra $P\subset pMp$ is {\it strongly non-amenable relative to $Q$} if for all non-zero projections $p_1\in P'\cap pMp,$ the von Neumann algebra $p_1P$ is non-amenable relative to $Q$.

For $B\subset M$ a von Neumann subalgebra, we have $L^2(M)\otimes_B L^2(M)\cong L^2(\langle M, e_B\rangle)$ as $M$-$M$-bimodules. Note that $B$ is amenable if and only if $_M L_2(M)\otimes_B L^2(M)_M \underset{weak}{\subset}$ $_M L^2(M)\otimes L^2(M)_M.$

Recall that a countable group $\Gamma$ is amenable if and only if every unitary representation of $\Gamma$ is weakly contained in the left regular representation (\cite[Theorem G.3.2]{BHV08}). The next lemma is the analogous statement for amenable von Neumann algebras. The result is likely well-known, but for a lack of reference, we include a proof.

\begin{lemma}\label{amenable}
Let $A$ be a tracial von Neumann algebra.  Then $A$ is amenable if and only if every $A$-$A$-bimodule $\mathcal K$ is weakly contained in the coarse bimodule $L^2(A)\otimes L^2(A).$ 

\end{lemma}

{\it Proof.} Suppose $A$ is amenable and let $\mathcal K$ be an $A$-$A$-bimodule. Then the trivial bimodule $_AL^2(A)_A$ is weakly contained in the coarse bimodule $_A L^2(A)\otimes L ^2(A)_A.$ Since $L^2(A)\otimes_A\mathcal K$ identifies with $\mathcal K$ as $A$-$A$ bimodules, we obtain that 

\begin{equation}\label{weak}
_A \mathcal K_A \underset{weak}{\subset} {}_A L^2(A)\otimes \mathcal K_A. 
\end{equation}

Now, since any right module of $A$ is contained in $\bigoplus_{\mathbb N} L^2(A)$ as a right $A$-submodule, we have that 
\begin{equation}\label{weak2}
_{\mathbb C}\mathcal K_{A} \underset{weak}{\subset} {} _{\mathbb C}L^2(A)_A. 
\end{equation}

Thus, \eqref{weak} and \eqref{weak2} implies that $_A\mathcal K_A$ is weakly contained in the coarse $A$-$A$-bimodule. The converse is clear by taking $\mathcal K=L^2(A)$, the trivial $A$-$A$-bimodule.
\hfill$\blacksquare$

We end this subsection by recording an immediate corollary of \cite[Lemma 2.6]{DHI16}. We provide a proof for the reader's convenience.

\begin{lemma}\cite[Lemma 2.6]{DHI16}\label{relative2}
Let $P$ and $Q$ be two von Neumann subalgebras of a tracial von Neumann algebra $(M,\tau)$. If $P$ is non-amenable relative to $Q$, then there exists a non-zero projection $z\in\mathcal N_{M}(P)'\cap M$ such that $Pz$ is strongly non-amenable relative to $Q$. 
\end{lemma}

{\it Proof.}
Using Zorn's lemma and a maximality argument, we can find a projection $z\in P'\cap M$ such that $Pz$ is strongly non-amenable relative to $Q$ and $P(1-z)$ is amenable relative to $Q$. Using \cite[Lemma 2.6]{DHI16} there exists $z_1\in\mathcal N_{M}(P)'\cap M$ such that $1-z\leq z_1$ and $Pz_1$ is amenable relative to $Q$. Therefore, $P(z_1-(1-z))$ is amenable relative to $Q$, which implies that $1-z=z_1\in \mathcal N_{M}(P)'\cap M$.

\hfill$\blacksquare$

\section{Intertwining of rigid algebras}\label{sectioni}

\subsection{The free product deformation for coinduced actions}\label{deformation}

In \cite{Io06a} Ioana introduced a malleable deformation for general Bernoulli actions, where the base is any tracial von Neumann algebra. This is a variant of the malleable deformation discovered by Popa \cite{Po03} in the case of Bernoulli actions with abelian or hyperfinite base and it was used in the context of coinduced actions in \cite{Dr15}.

Coinduced actions for tracial von Neumann algebras are defined as in Section \ref{coinduce}. More precisely, let $\Gamma$ be a countable group and let $\Sigma$ be a subgroup. Let $\phi:\Gamma/\Sigma \to \Gamma$ be a section. Define the cocycle $c:\Gamma\times\Gamma/\Sigma\to \Sigma$ by the formula $$c(g,i)=\phi^{-1}(gi)g\phi(i),$$ for all $g\in\Gamma$ and $i\in \Gamma/\Sigma.$\\
Let $\Sigma\overset{\sigma_{0}}{\curvearrowright} (A_0,\tau_0)$ be a trace preserving action, where $(A_0,\tau_0)$ is a tracial von Neumann algebra. We define an action $\Gamma\overset{\sigma}{\curvearrowright} A_0^{\Gamma/\Sigma}$, called the coinduced action of $\sigma_0$, as follows: 
$$\sigma_{g}((a_{i})_{i\in \Gamma/\Sigma})=(a'_{i})_{i\in \Gamma/\Sigma}, \text{  where } a'_{i}=(\sigma_0)_{c(g^{-1},i)^{-1}}(a_{g^{-1}i}).$$
Note that $\sigma$ is a trace preserving action of $\Gamma$ on the tracial von Neumann algebra $A_0^{\Gamma/\Sigma}$.

Consider the free product $A_0*L(\mathbb Z)$ with respect to the natural traces. Extend canonically $\sigma_0$ to an action on $A_0*L(\mathbb Z)$. Denote by $\tilde M=(A_0*L(\mathbb Z))^{\Gamma/\Sigma}\rtimes_{\sigma}\Gamma$ the corresponding crossed product of the coinduced action $\Gamma\overset{\sigma}{\curvearrowright} (A_0*L(\mathbb Z))^{\Gamma/\Sigma}$ of $\sigma_0.$

Take $u\in L(\mathbb{Z})$ the canonical generating Haar unitary. Let $h=h^{*}\in L(\mathbb{Z})$ be such that $u=$ exp$(ih)$ and set $u_{t}=$ exp$(ith)$ for all $t\in \mathbb{R}.$  Define the deformation $(\alpha_t)_{t\in\mathbb R}$ by automorphisms of $\tilde M$ by
\begin{center}
$\alpha_t(u_g)=u_g \, \,  \text{ and  } \, \alpha_{t}(\otimes_{h\in \Gamma/\Sigma}a_{h})=\otimes_{h\in \Gamma/\Sigma}$ Ad$(u_{t})(a_{h}),$\\
\end{center}
for all $g\in\Gamma, t\in\mathbb R$ and $\otimes_{h\in \Gamma/\Sigma}a_{h}\in (A_0*L(\mathbb Z))^{\Gamma/\Sigma}$ an elementary tensor.

\subsection{Spectral gap rigidity for coinduced actions}

\begin{theorem}\label{converge}
Let $\Gamma$ be an icc countable group and let $\Sigma$ be an almost malnormal subgroup. Let $\sigma_0$ be a pmp action of $\Sigma$ on a non-trivial standard probability space $(X_0,\mu_0)$. Denote by $M=L^\infty(X)\rtimes\Gamma$ the crossed-product von Neumann algebra of the coinduced action $\Gamma\overset{\sigma}\curvearrowright (X_0,\mu)^{\Gamma/\Sigma}$ associated to $\Sigma\overset{\sigma_0}\curvearrowright (X_0,\mu)$. Let $N$ be an arbitrary tracial von Neumann algebra and suppose $Q\subset p(M\bar\otimes N)p$ is a von Neumann subalgebra such that $Q'\cap p(M\bar\otimes N)p$ is strongly non-amenable relative to $1\otimes N$. 

Then, 
$$
\sup_{b\in\mathcal U(Q)}\|(\alpha_t\otimes \textrm{id})(b)-b\|_2 \text{ converges to 0 as } t\to 0.
$$
\end{theorem}

Theorem \ref{converge} and its proof are similar with other results from the literature \cite[Lemma 5.1]{Po06}, \cite[Corollary 4.3]{IPV10} and especially with \cite[Theorem 3.1]{BV12} (where the generalized Bernoulli action might have non amenable stabilizers) and with \cite[Theorem 2.6]{KV15} (which is another version of this result for coinduced actions).




{\bf Proof of Theorem \ref{converge}.}

Put $\mathcal M:=M\bar\otimes N$ and ${\mathcal {\tilde M}}:=\tilde M\bar\otimes N.$ The proof of this theorem goes along the same lines as the proof of \cite[Theorem 3.1]{BV12}. Therefore, instead of working with the bimodule $_{\mathcal M}{L^2(\tilde{\mathcal M}\ominus\mathcal M)}_{\mathcal M}$, we use the following $\mathcal M$-$\mathcal M$-submodule

\[
\mathcal K:=\overline{\text sp}
\left \{ (\otimes_{i\in\mathcal F}a_i) u_g\otimes n \middle |
 \begin{tabular}{ccc}
  $\mathcal F\subset \Gamma/\Sigma\ {\text  {with }} k\leq |\mathcal F|<\infty $, $n\in N$ and $g\in\Gamma$ \\
  $a_i\in A_0*L(\mathbb Z)$ for all $i\in\mathcal F$ \\
  $a_i\in (A_0*L(\mathbb Z))\ominus A_0$ for at least $k$ elements $i\in \mathcal F$
  \end{tabular}
\right \}.
\]

\vskip 0.05in
\noindent{\bf Claim 1.} The $\mathcal M$-$\mathcal M$-bimodule $\mathcal K$ is weakly contained in the bimodule $L^2(\mathcal M)\otimes_{1\otimes N}L^2(\mathcal M)$.

\vskip 0.05in
\noindent{\it Proof of Claim 1.}
Let $\mathcal A\subset A_0\ominus \mathbb C 1$ be an orthonormal basis of $L^2(A_0)\ominus \mathbb C 1$ and denote by $u$ the canonical Haar unitary of $L(\mathbb Z).$ Define the orthonormal set $\tilde A\subset L^2(A_0*L(\mathbb Z))\ominus L^2(A_0)$ by
\[\mathcal{\tilde A}:=\{ u^{n_1}a_1u^{n_2}a_2\dots u^{n_{k-1}}a_{k-1}u^{n_k}| k\ge 1, n_j\in\mathbb Z\setminus\{0\}, a_j\in\mathcal A \text{ for all } j \}
\]
This gives us the following orthogonal decomposition of $L^2(A_0*L(\mathbb Z))$ into $A_0$-$A_0$ submodules:
\begin{equation}\label{decomposition}
L^2(A_0*L(\mathbb Z))=L^2(A_0)\oplus\bigoplus_{a\in\mathcal{\tilde A}}\overline{A_0aA_0}.
\end{equation} 
If we denote 
$$
\mathcal C:=\{(\otimes_{i\in\mathcal F} c_i)\otimes 1| \mathcal F  \text{ finite } , k\leq|\mathcal F|<\infty, c_i\in\mathcal{\tilde A}, \text{ for all }i\in\mathcal F \},
$$
then the decomposition \eqref{decomposition} implies that the bimodule $\mathcal K$ can be written as the linear span  $\mathcal K=\overline{{\text {sp}}}_{c\in\mathcal C}\,\mathcal Mc\mathcal M$. To finish the proof of this claim, note that it is enough to consider an element $c \in\mathcal C$ and prove that the $ M$-$ M$-bimodule $\overline{\text {sp}}\, Mc M$ is weakly contained in the coarse bimodule $L^2( M)\otimes L^2( M).$

Let $c=(\otimes_{i\in\mathcal F} c_i)\otimes 1\in \mathcal C$. We denote by $\Gamma_0:=\{g\in\Gamma| gf=f, $ for all $f\in\mathcal F \}$, the stabilizer of $\mathcal F$ for the action $\Gamma\curvearrowright \Gamma/\Sigma$ and by $\Gamma_1:=\{g\in\Gamma|g\cdot\mathcal F=\mathcal F\}$, the normalizer of $\mathcal F$ for the same action. Since $\Sigma$ is $k$-almost malnormal and $\Gamma_0$ is a finite index subgroup of $\Gamma_1$, we obtain that $\Gamma_1$ is a finite group. \\
Denote $P=A\rtimes\Gamma_1$. Since $P$ is amenable, Lemma \ref{amenable} implies that the $P$-$P$-bimodule $\overline{\text {sp}}\, Mc M$ is weakly contained in the coarse bimodule $L^2(P)\otimes L^2(P).$
Thus, for each $\epsilon>0, F\subset\Gamma_1$ and $E\subset A$ finite subsets, there exist $\eta_1,\eta_2,...,\eta_n\in L^2(P)\otimes L^2(P)$ such that
\begin{equation}\label{partia0l}
|\langle au_g c (bu_h)^*,c \rangle-\sum_{i=1}^n \langle au_g\eta_i (bu_h)^*,\eta_i \rangle|\leq \epsilon,
\end{equation}
for all $g,h\in F$ and $a,b\in E.$

Using the canonical inclusion $L^2(P)\subset L^2( M)$, we obtain that $\langle au_g\eta_i (bu_h)^*,\eta_i \rangle$=0, for all $(g,h)\in(\Gamma\times\Gamma)\setminus (\Gamma_1\times\Gamma_1)$ and $a,b\in A.$ Note that also $\langle au_g c (bu_h)^*,c \rangle=0$, for all $(g,h)\in (\Gamma\times\Gamma)\setminus( \Gamma_1\times\Gamma_1)$ and $a,b\in A.$ Using these observations together with \eqref{partia0l}, we obtain that the $ M$-$ M$-bimodule $\overline{\text {sp}}\, Mc M$ is weakly contained in the coarse bimodule $L^2( M)\otimes L^2( M).$ This finishes the proof of the claim.
\hfill$\square$

Denote by $P_{\mathcal K}$ the orthogonal projection of $L^2(\mathcal{\tilde M})$ onto the closed subspace $\mathcal K$.

\vskip 0.05in
\noindent{\bf Claim 2.} $ \text{sup}_{b\in\mathcal U(Q)}\|P_{\mathcal K}((\alpha_t\otimes $id$)(b))\|_2$ converges to $0$ as $t\to 0.$

\vskip 0.05in
\noindent{\it Proof of Claim 2.}
Suppose the claim is false. Then there exist $\delta>0$, a sequence of positive numbers $t_n\to 0$, as $n\to\infty$, and a sequence of unitaries $b_n\in \mathcal U(Q)$ such that $\|P_{\mathcal K}((\alpha_{t_n}\otimes $id$)(b_n))\|_2\ge \delta$, for all $n\ge 1.$\\ 
Define $\xi_n=P_{\mathcal K}((\alpha_{t_n}\otimes $id$)(b_n))$. For all $x\in Q'\cap p\mathcal Mp$, we have $\|\xi_nx-x\xi_n\|\to 0$, as $n\to \infty.$ Note also that $\liminf_{n\to\infty}\|\xi_n\|_2\ge\delta$ and $\|x\xi_n\|_2\leq \|x\|_2$, for all $x\in\mathcal M.$ Then, \cite[Lemma 2.3]{Ho15} implies that there exists a projection $q\in \mathcal Z(Q'\cap p\mathcal Mp)$ such that the $\mathcal M$-$(Q'\cap p\mathcal Mp)q$ bimodule $L^2(\mathcal Mq)$ is weakly contained in $\mathcal K.$ Claim 1 implies now that the $\mathcal M$-$(Q'\cap p\mathcal Mp)q$ bimodule $L^2(\mathcal Mq)$ is weakly contained in the bimodule $L^2(\mathcal M)\otimes_{1\otimes N} L^2(\mathcal M).$ This implies that $(Q'\cap p\mathcal Mp)q$ is amenable relative to $1\otimes N$ inside $\mathcal M,$ which contradicts the hypothesis. This proves the claim.
\hfill$\square$

In order to finish the proof of the theorem we need a variant of Popa's transversality property. In the proof of \cite[Theorem 3.1]{BV12} it is proven the following fact for generalized Bernoulli actions: if $\text{sup}_{b\in\mathcal U(Q)}\|P_{\mathcal K}((\alpha_t\otimes $id$)(b))\|_2$ converges to $0$ as $t\to 0,$ then $\sup_{b\in\mathcal U(Q)}\|(\alpha_t\otimes \textrm{id})(b)-b\|_2 \text{ converges to 0 as } t\to 0.$ With the same proof we obtain the same result for coinduced actions. Claim 2 completes now the proof of the theorem.
\hfill$\blacksquare$

For $Q\subset M$ a von Neumann subalgebra, we define QN$_M(Q)\subset M$ to be the set of all elements $x\in M$ for which there exist $x_1,\dots,x_n,y_1,\dots,y_n$ satisfying $xQ\subset \sum_{i=1}^n Qx_i$ and $Qx\subset \sum_{i=1}^n y_i Q.$ The weak closure of QN$_M(Q)$ is called the {\it quasi-normalizer of $Q$ inside $M$} and note that it is a von Neumann subalgebra of $M$ which contains both $Q$ and $Q'\cap M.$ 

The proof of \cite[Theorem 4.2]{IPV10} carries over verbatim and gives us the following result.

\begin{theorem}\label{ipv}
Let $\Gamma$ be an icc countable group and let $\Sigma$ be an almost malnormal subgroup. Let $\sigma_0$ be a pmp action of $\Sigma$ on a non-trivial standard probability space $(X_0,\mu_0)$. Denote by $M=L^\infty(X)\rtimes\Gamma$ the crossed-product von Neumann algebra of the coinduced action $\Gamma\overset{\sigma}\curvearrowright (X_0,\mu)^{\Gamma/\Sigma}$ associated to $\Sigma\overset{\sigma_0}\curvearrowright (X_0,\mu)$. Let $N$ be a II$_1$ factor and suppose $Q\subset p(M\bar\otimes N)p$ is a von Neumann subalgebra. Denote by $P$ the quasi-normalizer of $Q$ in $p(M\bar\otimes N)p$.

If there exist $0<t<1$ and $\delta>0$ such that 
\[\tau(b^*(\alpha_t\otimes {\text{id}})(b))\ge\delta, \text{ for all } b\in\mathcal U(Q), 
\]
then one of the following statements is true:
\begin{itemize}
\item $Q\prec 1\otimes N$, 
\item $P\prec (A\rtimes\Sigma)\bar\otimes N$, 
\item there exists a unitary $u\in M\bar\otimes N$ such that $uPu^*\subset L(\Gamma)\bar\otimes N.$ 
\end{itemize}
\end{theorem}

\subsection{Controlling intertwiners and relative commutants}

In the Appendix of his PhD thesis \cite{Bo14}, Boutonnet has presented a unified approach to the notion of mixing for von Neumann algebras. As a consequence, we obtain results which give us good control over intertwiners between certain subalgebras of von Neumann algebras arising from coinduced actions.

\begin{definition}
Let $A\subset N\subset M$ be an inclusion of finite von Neumann algebras. We say that the inclusion $N\subset M$ is {\it mixing relative} to $A$ if for any sequence of unitaries $\{x_n\}\subset \mathcal U(N)$ with $\|E_A(yx_nz)\|_2\to 0$ for all $y,z\in N$, we have
$$
\|E_N(m_1x_nm_2)\|_2\to 0 \text{ for all } m_1,m_2\in M \ominus N. 
$$
\end{definition} 

\begin{proposition}\cite[Appendix A]{Bo14}\label{Mr}
Let $A\subset N\subset M$ be an inclusion of finite von Neumann algebras such that $N\subset M$ is mixing relative to $A$. 
Let $Q\subset pMp$ be a subalgebra such that $Q\nprec_M A$. Denote by $P$ the quasi-normalizer of $Q$ in $pMp.$ 
\begin{enumerate}
\item If $Q\subset N$, then $P\subset N$.
\item If $Q\prec N$, then there exists a non-zero partial isometry $v\in pM$ such that $vv^*\in P$ and $v^*Pv\subset N.$
\item If $N$ is a factor and if $Q\prec_{M}^s N$, then there exists a unitary $u\in \mathcal U(M)$ such that $uPu^*\subset N.$
\end{enumerate}

\end{proposition}

\begin{lemma}\label{Mc}
Let $\Sigma$ be a subgroup of a countable group $\Gamma$. Let $\Sigma\overset{\sigma_0}\car A_0$ be a tracial action on a non-trivial von Neumann algebra $A_0$ and let $\Gamma\overset{\sigma}\car A:=A_0^{\Gamma/\Sigma}$ be the coinduced action of $\sigma_0$. Let $\Gamma\car C$ be another tracial action. Then $C\rtimes\Gamma\subset (C\bar\otimes A)\rtimes\Gamma$ is mixing relative to $C\rtimes\Sigma.$
\end{lemma}

{\it Proof.}
Denote $\mathcal M :=(C\bar\otimes A)\rtimes\Gamma$ and $I:=\Gamma/\Sigma$. Let $\{x_n\}\subset\mathcal U(C\rtimes\Gamma)$ be a sequence of unitaries such that $\|E_{C\rtimes\Sigma}(yx_nz)\|_2\to 0,$ for all $y,z\in C\rtimes\Gamma.$ Let $a,b\in \mathcal M\ominus (C\rtimes\Gamma).$ We have to show that $\|E_{C\rtimes\Gamma}(ax_nb)\|_2\to 0.$ Since $E_{C\rtimes\Gamma}$ is $C\rtimes\Gamma$-bimodular, we can assume $a,b\in A.$ Moreover, we can suppose that there exist a finite subset $J\subset I$ and $j_0\in J$ such that $a,b=\otimes_{j\in J}b_j\in A_0^J$ with $b_{j_0}\in A_0\ominus \mathbb C.$ If $j_0=g_0\Sigma$ and $J=\{g_1\Sigma,\dots,g_n\Sigma\}$ note that $\Sigma_0:=\{g\in\Gamma| gj_0\in J\}=\cup_{i=1}^n g_i\Sigma g_0^{-1}$. Now, since $ax_nb=\sum_{g\in\Gamma} aE_{C}(x_nu_g^*)\sigma_g(b)u_g$, we have 
$$
\|E_{C\rtimes\Gamma}(ax_nb)\|_2^2=\sum_{g\in\Sigma_0}|\tau(a\sigma_{g}(b))|^2\|E_C(x_nu_g^*)\|_2^2\leq \|a\|_2^2\|b\|_2^2\sum_{i=1}^n\|E_{C\rtimes\Sigma}(u_{g_i}^*x_nu_{g_0})\|_2^2,
$$
which goes to zero because of the assumption. This proves the lemma.
\hfill$\blacksquare$

Proposition \ref{Mr} together with Lemma \ref{Mc} give the following result.

\begin{corollary}\label{unitary}
Let $\Sigma$ be a subgroup of a countable group $\Gamma$. Let $\Sigma\overset{\sigma_0}\car A_0$ be a tracial action on a non-trivial tracial von Neumann algebra $A_0$ and let $\Gamma\overset{\sigma}\car A:=A_0^{\Gamma/\Sigma}$ be the coinduced action of $\sigma_0$. Let $\Gamma\car C$ be another tracial action and let $N$ be an arbitrary factor. Define $\mathcal M:=(C\bar\otimes A)\rtimes\Gamma$. Suppose $Q\subset p(\mathcal M\bar\otimes N)p$ is a von Neumann subalgebra such that $Q\nprec (C\rtimes\Sigma)\bar\otimes N$. Denote by $P$ the quasi-normalizer of $Q$ inside $p(\mathcal M\bar\otimes N)p$.

\begin{enumerate}
\item If $Q\subset p((C\rtimes\Gamma)\bar\otimes N)p$, then $P\subset p((C\rtimes\Gamma)\bar\otimes N)p$.
\item If $Q\prec (C\rtimes\Gamma)\bar\otimes N$, then there exists a non-zero partial isometry $v\in p(\mathcal M\bar\otimes N)$ such that $vv^*\in P$ and $v^*Pv\subset (C\rtimes\Gamma)\bar\otimes N.$
\item If $Q\prec_{\mathcal M\bar\otimes N}^s (C\rtimes\Gamma)\bar\otimes N$, then there exists a unitary $u\in \mathcal U(\mathcal M\bar\otimes N)$ such that $uPu^*\subset (C\rtimes\Gamma)\bar\otimes N.$
\end{enumerate}

\end{corollary}

The proof of the following proposition is similar to \cite[Corollary 3.7]{Bo12a} and we leave it to the reader.

\begin{proposition}\label{remi}
Let $\Gamma\car C$ be a tracial action and denote $M_0=C\rtimes\Gamma.$ Let $\Sigma$ be an almost malnormal subgroup of $\Gamma$. Suppose $Q\subset pM_0p$ is a von Neumann subalgebra such that $Q\prec C\rtimes\Sigma$ and $Q\nprec C.$ Denote by $P$  the quasi-normalizer of $Q$ inside $pM_0p.$ \\
Then $P\prec C\rtimes\Sigma.$
\end{proposition}

\section{Rigidity coming from measure equivalence}\label{sectionr}

In this section we establish some results needed in the proof of Theorem \ref{2}. Throughout the section, we will work with coinduced actions satisfying the following:
\begin{assumption}\label{assumption} Let $\Sigma$ be a subgroup of a countable icc group $\Gamma$. Let $\sigma_0$ be a pmp action of $\Sigma$ on a non-trivial standard probability space $(X_0,\mu_0)$ and denote by $\sigma$ the coinduced action of $\Gamma$ on $X:=X_0^{\Gamma/\Sigma}$. Suppose:

\begin{itemize}
\item $\Gamma$ is a non-amenable icc group which is measure equivalent to a group $\Lambda_0$ for which the group von Neumann algebra $L(\Lambda_0)$ is not prime.
\item $\Sigma$ is almost malnormal.
\end{itemize}
\end{assumption}

Note that since $\Sigma$ is almost malnormal in $\Gamma$, we have that $[\Gamma:\Sigma]=\infty.$ Before stating the results of this section, we need to introduce some notation.
\begin{notation}\label{R}
The group von Neumann algebra $L(\Lambda_0)$ is not prime, therefore there exist von Neumann algebras $R_1$ and $R_2$, both not of type I, such that $L(\Lambda_0)=R_1\bar\otimes R_2$. Since $L(\Lambda_0)$ is diffuse and non-amenable, there exists $z_0\in\mathcal Z(L(\Lambda_0))$ such that $R_1z_0$ and $R_2z_0$ are diffuse and $L(\Lambda_0)z_0$ is non-amenable.

The group $\Gamma$ is measure equivalent to $\Lambda_0$. By \cite[Lemma 3.2]{Fu99}, $\Gamma$ and $\Lambda_0$ admit stably orbit equivalent free ergodic pmp actions. 
Thus, we may find a free ergodic pmp action $\Gamma\curvearrowright (Z_0,\nu)$ and $\ell\geq 1$, such that the following holds: consider the product action $\Gamma\times\mathbb Z/\ell\mathbb Z\curvearrowright (Z_0\times\mathbb Z/\ell\mathbb Z,\nu\times c)$, where $\mathbb Z/\ell\mathbb Z$ acts on itself by addition and $c$ denotes the counting measure on $\mathbb Z/\ell\mathbb Z$.  Then there exist a non-negligible measurable set $Y_0\subset Z_0\times\mathbb Z/\ell\mathbb Z$ and a free ergodic measure preserving action $\Lambda_0\curvearrowright Y_0$ such that $$\mathcal R(\Lambda_0\curvearrowright Y_0)=\mathcal R(\Gamma\times\mathbb Z/\ell\mathbb Z\curvearrowright Z_0\times\mathbb Z/\ell\mathbb Z)\cap (Y_0\times Y_0).$$

We put $C_0=L^{\infty}(Y_0), M_0=L^{\infty}(Z_0\times\mathbb Z/\ell\mathbb Z)\rtimes(\Gamma\times\mathbb Z/\ell\mathbb Z)$, $p=1_{Y_0}$, and note that
$C_0\rtimes\Lambda_0=pM_0p$. We identify $L^{\infty}(\mathbb Z/\ell\mathbb Z)\rtimes\mathbb Z/\ell\mathbb Z=\mathbb M_{\ell}(\mathbb C)$, and use this identification to write $M_0=C\rtimes\Gamma$, where $C=L^{\infty}(Z_0)\otimes\mathbb M_{\ell}(\mathbb C)$ and $\Gamma$ acts trivially on $\mathbb M_{\ell}(\mathbb C)$. 

Denote $A=L^\infty(X)$ and let $\{u_g\}_{g\in\Gamma}\subset (C\bar\otimes A)\rtimes\Gamma$ denote the canonical unitaries implementing the diagonal action of $\Gamma$ on $C\bar\otimes A.$
\end{notation}

\begin{remark} Throughout this section we will use many times the following easy observation (see \cite[Lemma 3.4]{Va08}). Let $P\subset pMp$ and $Q\subset qMq$ be von Neumann subalgebras of a tracial von Neumann algebra $(M,\tau).$ Then:
\begin{itemize}
\item if $p_0Pp_0\prec Q$ for a non-zero projection $p_0\in P$, then $P\prec Q.$
\item if $Pp'\prec Q$ for a non-zero projection $p'\in P'\cap pMp$, then $P\prec Q.$
\end{itemize}
\end{remark}

\begin{lemma}\label{d}
Let $w:\Gamma\to\mathcal U(A\bar\otimes N)$ be a cocycle for the action $\sigma\otimes$id, where $N$ is II$_1$ factor. Define the $*$-homomorphism $d:C\rtimes\Gamma\to (A\rtimes\Gamma)\bar\otimes N\bar\otimes (C\rtimes\Gamma)$ by $d(cu_g)=w_gu_g\otimes cu_g, g\in\Gamma, c\in C.$ Let $Q\subset pM_0p$ be a subalgebra and let $\Sigma_0\subset \Gamma$ be a subgroup. The following hold:
\begin{enumerate}
\item If $Q\nprec C$ , then $d(Q)\nprec 1\otimes N\bar\otimes (C\rtimes\Gamma).$
\item If $[\Gamma:\Sigma_0]=\infty$, then $d(L(\Lambda_0))\nprec (A\rtimes\Sigma_0)\bar\otimes N\bar\otimes (C\rtimes\Gamma).$
\item If $Q$ is non-amenable, then $d(Q)$ is non-amenable relative to $1\otimes N\bar\otimes (C\rtimes\Gamma).$
\end{enumerate}
\end{lemma}

{\it Proof.} Denote $\mathcal M=(A\rtimes\Gamma)\bar\otimes N\bar\otimes M_0$ and $\mathcal N=1\otimes N\bar\otimes M_0.$\\ 
(1) Let $\{u_n\}_{n\ge 1}\subset \mathcal U(Q)$ be a sequence of unitaries such that $\|E_C(u_nu_g)\|_2\to 0$, for all $g\in\Gamma$. We claim that 
$$
\|E_{1\otimes N\bar\otimes M_0}(xd(u_n)y)\|_2\to 0, \text{ for all } x,y\in\mathcal M.
$$
Since $E_{\mathcal N}$ is $\mathcal N$-bimodular, by Kaplansky's density theorem we may assume $x=au_g\otimes 1\otimes 1$, $y=bu_h\otimes 1\otimes 1$ for some $a,b\in A$ and $g,h\in\Gamma$. Then for all $n\ge 1,$ we have
\[xd(u_n)y=\sum_{k\in\Gamma} a\sigma_g(w_{k})\sigma_{gk}(b)u_{gkh}\otimes E_C(u_nu_{k}^*)u_k.
\] 
Therefore, $\|E_{\mathcal N}(xd(u_n)y)\|_2\leq \|a\|\|b\|\|E_C(u_nu^*_{g^{-1}h^{-1}})\|_2\to 0.$

(2)  Assume $d(L(\Lambda_0))\prec (A\rtimes\Sigma_0)\bar\otimes N\bar\otimes M_0$. Since $d(C_0)\subset 1\otimes 1\otimes C_0$, we obtain $d(C_0\rtimes\Lambda_0)\prec (A\rtimes\Sigma_0)\bar\otimes N\bar\otimes M_0$. Therefore $d(L(\Gamma))\prec (A\rtimes\Sigma_0)\bar\otimes N\bar\otimes M_0,$ which implies $L(\Gamma)\prec L(\Sigma_0)$. Indeed, suppose by contrary that $L(\Gamma)\nprec L(\Sigma_0)$. Then there exists a sequence $u_n\in\mathcal U(L(\Gamma))$ such that $\|E_{L(\Sigma_0)}(xv_ny)\|_2\to 0$, for all $x,y\in L(\Gamma).$ We would like to prove that
\begin{equation}\label{ex}
\|E_{(A\rtimes\Sigma_0)\bar\otimes N\bar\otimes M_0}(xd(u_n)y)\|_2\to 0,
\end{equation}
for all $x,y\in (A\rtimes\Gamma)\bar\otimes N\bar\otimes M_0.$ For proving \eqref{ex}, it is enough to consider $x=u_g\otimes 1\otimes 1$ and $y=u_h\otimes 1\otimes 1,$ with $g,h\in\Gamma.$ In this case one can check that 
\[ \|E_{(A\rtimes\Sigma_0)\bar\otimes N\bar\otimes M_0}(xd(u_n)y)\|_2=\|E_{L(\Sigma_0)}(u_gu_nu_h)\|_2,
\]
which goes to $0$. Therefore \eqref{ex} is proven and we obtain that $d(L(\Gamma))\nprec (A\rtimes\Sigma_0)\bar\otimes N\bar\otimes M_0,$ contradiction. 

Thus $L(\Gamma)\prec L(\Sigma_0)$, which implies that $\Sigma_0$ has finite index in $\Gamma$ by \cite[Lemma 2.5]{DHI16}. 

(3) Suppose by contrary that $d(Q)$ is amenable relative to $\mathcal N.$ Then there exists a positive linear functional $\phi:d(p)\langle \mathcal M,e_{\mathcal N}\rangle d(p)\to \mathbb C$ such that $\phi_{|d(p)\mathcal Md(p)}=\tau$ and $\phi$ is $d(Q)$-central. Define now 
$\varphi:p\langle M_0,e_{\mathbb C} \rangle p\to \mathbb C$ by 
$$\varphi (\sum_{i=1}^N m_i e_{\mathbb C}n_i)=\phi(\sum_{i=1}^N d(m_i)e_{\mathcal N}d(n_i)),$$
where $N\ge 1$, $m_i,n_i\in  M_0, i\in\{1,\dots,N\}$. Note that $\varphi$ is a well defined positive linear functional. Indeed, suppose $\sum_{i=1}^N m_ie_{\mathbb C}n_i=0$, with $m_i,n_i\in M_0$, for all $1\leq i\leq N.$ This implies $\sum_{i=1}^N d(m_i)\tau(n_i)=0$.
Since $E_{\mathcal N}(d(m))=\tau(m)$, for all $m\in M_0$, we obtain $\sum_{i=1}^N d(m_i) E_{\mathcal N}(d(n_i))=0,$ which implies $\sum_{i=1}^N d(m_i) e_{\mathcal N}d(n_i)=0.$ Therefore, $\varphi$ is a positive linear functional which is $Q$-central and $\varphi_{|pM_0p}=\tau$. We obtain $Q$ is amenable, contradiction.
\hfill$\blacksquare$

Denote by $\mathcal U_{fin}$ the class of Polish groups which arise as closed subgroups of the unitary groups of II$_1$ factors \cite{Po05}. In particular, all countable discrete groups and all compact Polish groups belong to $\mathcal U_{fin}.$

\begin{theorem}{(Cocycle superrigidity.)}\label{cocycle}
Let $\Gamma\curvearrowright X$ be as in Assumption \ref{assumption}. Then any cocycle $w:\Gamma\times X\to \Lambda$ valued in a group $\Lambda\in \mathcal U_{fin}$ untwists, i.e. there exists a measurable map $\varphi:X\to \Lambda$ and a group homomorphism $d:\Gamma\to\Lambda$ such that $w(g,x)=\varphi(gx)d(g)\varphi(x)^{-1}$ for all $g\in\Gamma$ and a.e. $x\in X$. 
\end{theorem}
This result was proven in \cite{PS09} for Bernoulli actions using deformations obtained from closable derivations. In our case, we will provide a direct proof for Theorem \ref{cocycle} which uses only the free product deformation $\alpha_t$ defined in Section \ref{deformation}.

{\it Proof.}
Define $A:=L^\infty(X)$ and let $N$ be a II$_1$ factor such that $\Lambda\subset \mathcal U(N).$ We associate to $w:\Gamma\times X\to\mathcal U( N)$ the cocycle $w:\Gamma\to\mathcal U(A\bar\otimes N)$, given by $w_g(x)=w(g,g^{-1}x).$ Define $Q=\{w_gu_g\}_{g\in\Gamma}''.$

{\bf Claim.} We have $$\sup_{b\in\mathcal U(Q)}\|(\alpha_t\otimes \textrm{id})(b)-b\|_2 \text{ converges to 0 as } t\to 0.
$$

{\it Proof of the Claim.}
As in Lemma \ref{d} we define the $*$-homomorphism $d:C\rtimes\Gamma\to (A\rtimes\Gamma)\bar\otimes N\bar\otimes (C\rtimes\Gamma)$ by $d(cu_g)=w_gu_g\otimes cu_g, g\in\Gamma, c\in C.$ Denote $\mathcal M=(A\rtimes\Gamma)\bar\otimes N\bar\otimes M_0.$ 
%
Without loss of generality assume that $R_1z_0$ is non-amenable. Lemma \ref{d} implies that  $d(R_1z_0)$ is non-amenable relative to $1\otimes N\bar\otimes (C\rtimes\Gamma)$.
By Lemma \ref{relative2} there exists a non-zero projection $z\in\mathcal N_{d(z_0)\mathcal M d(z_0)} d(R_1z_0)'\cap d(z_0)\mathcal M d(z_0)$ such that $d(R_1)z$ is strongly non-amenable relative to $1\otimes N\bar\otimes (C\rtimes\Gamma)$. Using Theorem \ref{converge} we obtain that 
$$\sup_{b\in\mathcal U(d(R_2)z)}\|(\alpha_t\otimes \textrm{id}\otimes \textrm{id})(b)-b\|_2 \text{ converges to 0 as } t\to 0.
$$
and therefore by Theorem \ref{ipv} we obtain that one of the following hold:
\begin{enumerate}
\item $d(R_2)z\prec 1\otimes N\bar\otimes (C\rtimes\Gamma)$,
\item $d(L(\Lambda_0))z\prec (A\rtimes\Sigma)\bar\otimes N\bar\otimes (C\rtimes\Gamma),$ 
\item $d(L(\Lambda_0))z\prec L(\Gamma)\bar\otimes N\bar\otimes (C\rtimes\Gamma).$ 
\end{enumerate}
Note that (1) and (2) are not possible by Lemma \ref{d} since $R_2z_0$ is diffuse and  $[\Gamma:\Sigma]=\infty.$
Therefore (3) is true. 

Now, together with the remark that $d(C)\subset 1\otimes 1\otimes C$ we obtain that
$d(C\rtimes\Gamma)\prec L(\Gamma)\bar\otimes N\bar\otimes (C\rtimes\Gamma)$. One can check directly this fact or use \cite[Lemma 2.3]{BV12}. Proceeding in the same way, we obtain actually $d(C\rtimes\Gamma)\prec_{\mathcal M}^s L(\Gamma)\bar\otimes N\bar\otimes (C\rtimes\Gamma)$. Lemma \ref{d} implies that $d(C\rtimes\Gamma)\nprec L(\Sigma)\bar\otimes N\bar\otimes (C\rtimes\Gamma)$, so by Corollary \ref{unitary} we obtain
$$\sup_{b\in\mathcal U(Q)}\|(\alpha_t\otimes \textrm{id})(b)-b\|_2 \text{ converges to 0 as } t\to 0.
$$
\hfill$\square$

Using a result which goes back to Popa \cite{Po05}, the claim implies that the cocycle $w$ untwists (see \cite[Theorem 2.15]{Dr15}, the proof of \cite[Proposition 3.2]{Dr15} and \cite[Remark 3.3]{Dr15}).
\hfill$\blacksquare$

\begin{theorem}\label{delta}
Let $\Gamma\curvearrowright X$ be as in Assumption \ref{assumption} and supppose that $\Sigma$ is amenable. Let $\Lambda\curvearrowright B$ be a tracial action on a non-trivial abelian von Neumann algebra $B$ such that $A\rtimes\Gamma=B\rtimes\Lambda.$ Denote by $\Delta:B\rtimes\Lambda\to (B\rtimes\Lambda)\bar\otimes L(\Lambda)$ the comultiplication $\Delta(bv_\lambda)=bv_\lambda\otimes v_\lambda$ for all $b\in B$ and $\lambda\in \Lambda$ (we let $\{v_\lambda\}_{\lambda\in\Lambda}\subset B\rtimes\Lambda$ denote the canonical unitaries implementing the action of $\Lambda$ on $B$).

Then there exists a unitary $u\in\mathcal U((A\rtimes\Gamma)\bar\otimes (A\rtimes\Gamma))$ such that $$u\Delta(L(\Gamma))u^*\subset L(\Gamma\times\Gamma).$$
\end{theorem}

Define $M:=(C\bar\otimes A)\rtimes\Gamma$ and $\theta:M\to M\bar\otimes M\bar\otimes M$ by $\theta(cau_g)=cu_g\otimes\Delta(au_g),$ for all $c\in C, a\in A$ and $g\in\Gamma.$ In the following lemma we record some properties of the unital $*$-homomorphism $\theta$ which are similar to the ones of \cite[Lemma 10.2]{Io10}.

\begin{lemma}\label{diagonal3}

Let $Q\subset qMq$. The following hold:
\begin{enumerate}
\item If $Q$ is diffuse, then $\theta(Q)\nprec M\bar\otimes 1\bar\otimes M$.
\item If $Q\nprec B$, then $\theta(Q)\nprec M\bar\otimes M\otimes 1.$
\item If $Q$ has no amenable direct summand, then $\theta(Q)$ is strongly non-amenable relative to $M\bar\otimes M\otimes 1$ and $M\bar\otimes 1\bar\otimes M.$
\end{enumerate}
\end{lemma}

We continue now with the proof of Theorem \ref{delta} and we will give the proof of Lemma \ref{diagonal3} at the end of this section. 

{\it Proof of Theorem \ref{delta}.} 
Without loss of generality we can assume that $R_1z_0$ is non-amenable. Take $z\in \mathcal Z(R_1z_0)$ such that $R_1z$ has no amenable direct summand.

{\bf Claim 1.} We have $\sup_{b\in\mathcal U(\Delta(L(\Gamma))}\|( \textrm{id}\otimes \alpha_t)(b)-b\|_2 \text{ converges to 0 as } t\to 0.$

{\it Proof of Claim 1.}
Note that $\theta (R_1z)$ is strongly non-amenable relative to $M\bar\otimes M\otimes 1$ by Lemma \ref{diagonal3}. Therefore by Theorem \ref{converge}, we obtain
\begin{equation}\label{conv}
\sup_{b\in\mathcal U(\theta (R_2z))}\|(\textrm{id}\otimes \textrm{id}\otimes \alpha_t)(b)-b\|_2 \text{ converges to 0 as } t\to 0.
\end{equation}
 
Using Theorem \ref{ipv} we obtain that one of the following three conditions holds:
\begin{enumerate}
\item $\theta (R_2z)\prec M\bar\otimes M\otimes 1$,
\item $\theta (L(\Lambda_0)z)\prec M\bar\otimes M\bar\otimes (A\rtimes\Sigma)$,
\item there exists a unitary $u\in \mathcal M$ such that $u\theta (L(\Lambda_0)z)u^*\subset M\bar\otimes M\bar\otimes L(\Gamma).$
\end{enumerate}
  
If (1) holds, Lemma \ref{diagonal3} implies $R_2z \prec_M B$. By applying \cite[Lemma 3.5]{Va08}, we obtain $B\prec_M zMz\cap (R_2z)'$. Note that if $R_2z\prec_M C\rtimes\Sigma$, Proposition \ref{remi} implies that $L(\Lambda_0)\prec C\rtimes\Sigma$. Using \cite[Lemma 2.3]{BV12} we deduce that $C\rtimes\Gamma\prec C\rtimes\Sigma$. This is a contradiction since $[\Gamma:\Sigma]=\infty.$\\
Therefore $R_2z\nprec_M C\rtimes\Sigma$ and Corollary \ref{unitary} implies that $zMz\cap (R_2z)'\subset C\rtimes\Gamma,$ so $B\prec_M C\rtimes\Gamma.$ On the other hand, since $B\subset A\rtimes\Gamma,$ we obtain $B\prec_{A\rtimes\Gamma} L(\Gamma).$ Proposition \ref{remi} implies that $B\nprec_{A\rtimes\Gamma} L(\Sigma)$. Finally, using Corollary \ref{unitary} we obtain that $A\rtimes\Gamma\prec_{A\rtimes\Gamma} L(\Gamma)$, which is a contradiction.

Now, if (2) holds, we obtain $\theta (L(\Lambda_0))\prec M\bar\otimes M\bar\otimes (A\rtimes\Sigma)$. Together with $\theta(C)\subset C\otimes 1\otimes 1$, we obtain $\theta(M_0)\prec M\bar\otimes M\bar\otimes (A\rtimes\Sigma).$ Since $\Sigma$ is amenable, it implies that $\theta(M_0)$ is not strongly non-amenable relative to $M\bar\otimes M\otimes 1.$ Now, $M_0$ is a factor, so Lemma \ref{diagonal3} gives that $M_0$ is amenable, which is a contradiction. 

Thus, (3) holds. Since $\theta(C_0)\subset C_0\otimes 1\otimes 1$, we obtain  
\[\theta(M_0)\prec M\bar\otimes M\bar\otimes L(\Gamma)
\] 
With the same computation, we obtain $\theta(M_0)\prec_{M\bar\otimes M\bar\otimes M}^s M\bar\otimes M\bar\otimes L(\Gamma).$

Lemma \ref{diagonal3} implies that $\theta(M_0)\nprec M\bar\otimes M\bar\otimes L(\Sigma)$ since $\Sigma$ is amenable and $M_0$ is a factor. By Corollary \ref{unitary} we obtain that $\sup_{b\in\mathcal U(\Delta(L(\Gamma))}\|( \textrm{id}\otimes \alpha_t)(b)-b\|_2 \text{ converges to 0 as } t\to 0.$ 
\hfill$\square$

{\bf Claim 2.} 
We have $\sup_{b\in\mathcal U(\Delta(L(\Gamma)))}\|(\alpha_t\otimes\textrm{id})(b)-b\|_2 \text{ converges to 0 as } t\to 0.$ 

{\it Proof of Claim 2.}
As in Claim 1, by applying Lemma \ref{diagonal3}, Theorem \ref{converge} and Theorem \ref{ipv} we obtain that one of the following conditions hold: 

\begin{enumerate}
\item $\theta (R_2z)\prec M\bar\otimes 1\bar\otimes M$,
\item $\theta (L(\Lambda_0)z)\prec M\bar\otimes (A\rtimes\Sigma)\bar\otimes M$,
\item there exists a unitary $u\in \mathcal M$ such that $u\theta (L(\Lambda_0)z)u^*\subset M\bar\otimes L(\Gamma)\bar\otimes M.$
\end{enumerate}

Note that by Lemma \ref{diagonal3}, (1) is not possible since $R_2z$ is diffuse. As before, (2) is not possible, which implies (3) holds true and by reasoning as before we obtain the claim.
\hfill$\square$

Notice that $\Delta(L(\Gamma))$ is a factor since $\Gamma$ is icc. Using Claim 1 and 2 and by applying twice Theorem \ref{ipv} and \cite[Lemma 10.2]{IPV10} we obtain the conclusion.
\hfill$\blacksquare$

{\it Proof of Lemma \ref{diagonal3}.} The proofs of (1) and (2) are similar to the proof of Lemma \ref{d}.1 (see also the proof of \cite[Lemma 10.2]{IPV10}). For proving (3), denote $\mathcal M:=M\bar\otimes M\bar\otimes (A\rtimes\Gamma)$ and $\psi:M\to M\bar\otimes M$, by $\psi(cau_g)=cu_g\otimes au_g$ for all $c\in C, a\in A$ and $g\in\Gamma.$

\vskip 0.05in
\noindent{\bf Claim 1.} We have $_{\mathcal M}L^2(\mathcal M)\otimes_{M\bar\otimes M\otimes 1} L^2(\mathcal M)_{\theta( M)}\underset{weak}{\subset} {}  {_{\mathcal M}}L^2(\mathcal M)\otimes L^2(A\rtimes\Gamma)_{\psi(M)_{1,4}}$.\\
 (here we consider that $\psi(M)\subset M\bar\otimes M$ acts to the right on $L^2(M)\otimes L^2(M)\otimes L^2(M)\otimes L^2(M)$ on the first and fourth positions.)

\vskip 0.05in
\noindent{\it Proof of the Claim 1.}
Note that we have the identification
$$_{\mathcal M}L^2(\mathcal M)\otimes_{M\bar\otimes M\otimes 1} L^2(\mathcal M)_{\theta( M)}\simeq  _{\mathcal M_{1,2,3}}L^2(M\bar\otimes M\bar\otimes (A\rtimes\Gamma)\bar\otimes (A\rtimes\Gamma))_{\theta( M)_{1,2,4}}
$$
as $\mathcal M$-$M$-bimodules. 
Therefore, it is enough to show that 
$$_{M\bar\otimes M\otimes 1}L^2(M\bar\otimes M\bar\otimes (A\rtimes\Gamma))_{\theta(M)} \underset{weak}{\subset} {} _{M\bar\otimes M\otimes 1}L^2(M\bar\otimes M\bar\otimes (A\rtimes\Gamma))_{\psi(M)_{1,3}}.
$$

Let $\mathcal B$ be an orthonormal basis for $L^2(B)$ and note that we have the following orthogonal decomposition into $(M\bar\otimes M)$-$M$-bimodules:
\[L^2(M\bar\otimes M\bar\otimes (A\rtimes\Gamma))=\bigoplus_{b\in\mathcal B}\overline {\text {sp}}\, (M\bar\otimes M\otimes 1 )(1\otimes 1\otimes b)\, \theta (M)
\]

First, notice that for a fixed $b\in\mathcal B$ we have  
$${\overline{\text {sp}}}\, (M\bar\otimes M\otimes 1 )(1\otimes 1\otimes b)\, \theta (M)\simeq _{(M\bar\otimes M)_{1,2}}{L^2(M)\otimes L^2(M)\otimes_{B} L^2(A\rtimes\Gamma)}_{\psi(M)_{1,3}}
$$
as $(M\bar\otimes M)$-$M$-bimodules. 
Indeed, let $m_1,m_2,m_3\in M$ and let us prove that
\begin{equation}\label{bimodule}
\langle (m_1\otimes m_2\otimes 1) (1\otimes 1\otimes b)\theta(m_3),1\otimes 1\otimes b\rangle=
\langle (m_1\otimes m_2\otimes 1) (1\otimes 1\otimes_B 1) \psi(m_3),1\otimes1\otimes_B 1\rangle
\end{equation}
We may assume $m_3=cau_g$ for some $c\in C, a\in A$ and $g\in\Gamma$. Write $au_g=\sum_{l\in\Lambda} b_lv_l\in B\rtimes\Lambda$, with $b_l\in B$ for all $l\in\Lambda.$ Therefore, the LHS of \eqref{bimodule} equals to
$$
\tau((m_1\otimes m_2\otimes b^*b) \theta(m_3))=\tau(m_1cu_g\otimes ((m_2\otimes b^*b)\Delta(au_g)))=
\tau(m_1cu_g)\tau(m_2b_e).
$$
On the other hand, the RHS of \eqref{bimodule} equals to
$$
\tau((m_1\otimes E_B (m_2))\psi(m_3))=\tau (m_1cu_g\otimes E_B(m_2)au_g)=\tau(m_1cu_g)\tau(m_2b_e),
$$
which proves \eqref{bimodule}.

Now since $B$ is amenable, we obtain that 
$$_{(M\bar\otimes M)_{1,2}}{L^2(M)\otimes L^2(M)\otimes_{B} L^2(A\rtimes\Gamma)}_{\psi(M)_{1,3}}\underset{weak}{\subset} {} _{(M\bar\otimes M)_{1,2}}{L^2(M)\otimes L^2(M)\otimes L^2(A\rtimes\Gamma)}_{\psi(M)_{1,3}}.
$$
This finishes the proof of the claim.    
\hfill$\square$

\vskip 0.05in
\noindent{\bf Claim 2.} We have $ {_{\mathcal M}}L^2(\mathcal M)\otimes L^2(M)_{\psi(M)_{1,4}}\underset{weak}{\subset} {} _\mathcal M L^2(\mathcal M)\otimes L^2(M)_M$.

\vskip 0.05in
\noindent{\it Proof of the Claim 2.}
First, note that it is enough to prove  
$$_{M}{L^2(M)\otimes L^2(M)}_{\psi(M)}\underset{weak}{\subset} {} _M L^2(M)\otimes L^2(M)_M. 
$$
Let $\mathcal C$ be an orthonormal basis for $L^2(C)$ and note that we have the following orthogonal decomposition into $M$-$M$-bimodules:
\[
L^2(M)\otimes L^2(M)=\bigoplus_{c\in\mathcal C}\overline {\text {sp}}\, M\,(1\otimes c)\,d(M).
\]
Note that $ \overline {\text {sp}}\, M\,(1\otimes c)\,d(M)    \cong    L^2(M)\otimes_C L^2(M)$ as $M$-$M$-bimodules.
Indeed, let us take $m_1=c_1a_1u_{g_1}$, $m_2=c_2a_2u_{g_2}$, and note that 
\[
\langle m_1 (1\otimes c)\psi (m_2), 1\otimes c \rangle=\langle c_1a_1u_{g_1}c_2u_{g_2}\otimes ca_2u_{g_2},1\otimes c \rangle=\delta_{g_1,e}\delta_{g_2,e}\tau(c_1c_2)\tau(a_1)\tau(a_2)
\]
and
\[
\langle m_1e_Cm_2,e_C \rangle=\tau(E_C(c_1a_1u_{g_1})c_2a_2u_{g_2})=\delta_{g_1,e}\delta_{g_2,e}\tau(c_1c_2)\tau(a_1)\tau(a_2).
\]
This implies that $ \overline {\text {sp}}\, M\,(1\otimes c)\,\psi(M)   \cong    L^2(M)\otimes_C L^2(M)$ as $M$-$M$-bimodules. Since $C$ is amenable, the claim is proven.  
\hfill$\square$

Now, assume that  $\theta(Q)$ is not strongly non-amenable relative to $M\bar\otimes M\otimes 1.$ Then there exists a non-zero projection $p\in \theta(Q)'\cap \theta(q)\mathcal M\theta(q)$ such that 
$$_{\mathcal M}{L^2(\mathcal M p)}_{\theta(Q)}  \underset{weak}{\subset} {}  _{\mathcal M}L^2(\mathcal M)\otimes_{M\bar\otimes M\otimes 1} L^2(\mathcal M)_{\theta( Q)}.
$$

Using Claim 1 and 2, we obtain now that $_{\mathcal M}{L^2(\mathcal M p)}_{\theta(Q)}  \underset{weak}{\subset} {}  _\mathcal M L^2(\mathcal M)\otimes L^2(Q)_Q.$

Take $z\in Q$ such that $\theta (z)$ is the support projection of $E_{\theta(Q)}(p).$ Note that $z$ is a non-zero central projection in $Q$ and that $\theta$ embeds the trivial $Qz$-$Qz$-bimodule into $_{\theta(Qz)}{L^2({\theta({Qz})})}_{\theta(Qz)}$. Therefore, $_{Qz} L^2(Qz)_{Qz} \underset{weak}{\subset} {}   _{\theta(Qz)} L^2(\mathcal M)\otimes L^2(Qz)_{Qz}$.  Finally, we obtain $_{Qz} L^2(Qz)_{Qz} \underset{weak}{\subset} {}   _{Qz} L^2(Qz)\otimes L^2(Qz)_{Qz}$, which means that $Qz$ is amenable, contradiction. 

In a similar way, one can prove that $\theta(Q)$ is strongly non-amenable relative to $M\bar\otimes 1\bar\otimes M.$ This ends the proof.
\hfill$\blacksquare$

\section{Intertwining of abelian subalgebras}\label{sectionio}
Throughout this section we will use the following notation. Let $\Gamma$ be a countable group. Let $\Sigma$ be an almost malnormal subgroup and let $\sigma_0$ be a tracial action of $\Sigma$ on a non-trivial abelian von Neumann algebra $A_0$. Denote by $\sigma$ the coinduced action of $\Gamma$ on $A:=A_0^{\Gamma/\Sigma}$. Finally, denote $M=A\rtimes\Gamma.$

The next result is a localization theorem for coinduced actions which goes back to\cite[Theorem 6.1]{Io10}. The form presented in this paper is very similar to \cite[Theorem 5.1]{IPV10}, but written with coinduced actions instead of generalized Bernoulli ones.

\begin{theorem}\label{clustering}
Assume that $D\subset M\bar\otimes M$ is an abelian von Neumann subalgebra which is normalized by a group of unitaries $(\gamma(s))_{s\in\Lambda}$ that belong to $L(\Gamma)\bar\otimes L(\Gamma)$. Denote by $P$ the quasi-normalizer of $D$ inside $M\bar\otimes M.$ We make the following assumptions:

\begin{enumerate}
\item $D\nprec M\otimes 1$ and $D\nprec 1\otimes M$,
\item $P\nprec M\bar\otimes (A\rtimes\Sigma)$ and $P\nprec (A\rtimes\Sigma)\bar\otimes M$,
\item $P\nprec M\bar\otimes L(\Gamma)$ and $P\nprec L(\Gamma)\bar\otimes M$,
\item $\gamma(\Lambda)''\nprec L(\Gamma)\bar\otimes L(\Sigma)$ and $\gamma(\Lambda)''\nprec L(\Sigma)\bar\otimes L(\Gamma)$.
\end{enumerate}

Define $C:=D'\cap (M\bar\otimes M).$ Then for every non-zero projection $q\in\mathcal Z(C)$ we have $Cq\prec A\bar\otimes A.$
\end{theorem}

The proof is identically with the one of \cite[Theorem 5.1]{IPV10}, since essentially the same computations still hold once we replace generalized Bernoulli actions by coinduced ones.

Next, we obtain a similar statement if one considers an abelian von Neumann algebra in $M$ and not in $M\bar\otimes M$.

\begin{theorem}\label{clustering2}
Assume that $D\subset M$ is an abelian von Neumann subalgebra which is normalized by a group of unitaries $(\gamma(s))_{s\in\Lambda}$ that belong to $ L(\Gamma)$. Denote by $P$ the quasi-normalizer of $D$ inside $M.$ We make the following assumptions:

\begin{enumerate}
\item $D$ is diffuse,
\item $P\nprec A\rtimes\Sigma$,
\item $P\nprec L(\Gamma) $,
\item $\gamma(\Lambda)''\nprec L(\Sigma)$.
\end{enumerate}

Define $C:=D'\cap M.$ Then for every non-zero projection $q\in\mathcal Z(C)$ we have $Cq \prec A.$
\end{theorem}

As noticed in \cite{Io10}, we obtain as a corollary a weaker version of Popa's conjugacy criterion adapted in this case to coinduced actions.

\begin{theorem}\label{criterion}
Suppose $\Gamma$ is icc and $\Sigma$ is amenable. Let $\Lambda\curvearrowright B$ be another tracial action of a countable group $\Lambda$ on a non-trivial abelian von Neumann algebra $B$ such that $M=A\rtimes\Gamma=B\rtimes\Lambda$ and $L(\Lambda)\subset L(\Gamma)$.

Then $B\prec A.$
\end{theorem} 

{\it Proof.}
The proof is a direct application of Theorem \ref{clustering2}. Note that the quasi-normalizer of the abelian algebra $B$ is $M$. Now, notice that if $M\prec A\rtimes\Sigma,$ by \cite[Lemma 2.5.1]{DHI16} we obtain that $[\Gamma:\Sigma]<\infty$. This is not possible since $\Sigma$ is almost malnormal in $\Gamma$. Also $L(\Lambda)\nprec L(\Sigma)$ since $\Sigma$ is amenable and therefore we obtain $B\prec A.$

\hfill$\blacksquare$

\section{Proof of the main results}\label{sectionp}

In \cite{Io10}, Ioana has proven that any Bernoulli action of an arbitrary icc property (T) group is W$^*$-superrigid. The strategy of his proof was successfully applied also in \cite{IPV10} and \cite{Bo12}. 

\subsection{\bf A general method for obtaining W$^*$-superrigidity.}
Using Ioana's proof, we identify a couple of steps for proving that a certain free ergodic pmp action $\Gamma\curvearrowright X$ is W$^*$- superrigid (see also the introduction of \cite{Bo12}). Consider an arbitrary free ergodic pmp action $\Lambda\curvearrowright Y$ such that $M:=A\rtimes\Gamma=B\rtimes\Lambda,$ where $A=L^\infty(X)$ and $B=L^\infty(Y).$
 Define the comultiplication $\Delta:M\to M\bar\otimes L(\Lambda)$ by $\Delta(bv_\lambda)=bv_\lambda\otimes v_\lambda,$ for all $ b\in B,\lambda\in\Lambda$, where we denote by $v_\lambda,\lambda\in\Lambda$, the canonical unitaries corresponding to the action of $\Lambda.$

\begin{enumerate}[{\bf Step 1.}]
\item One has to show that $\Gamma\curvearrowright X$ is OE superrigid. From now on, using Singer's result \cite{Si55}, it is enough to assume that $B$ is not unitarily conjugated to $A$ in $M$, which is equivalent to $B\nprec_M A$ \cite[Theorem A.1]{Po06b}.

\item One can also assume that there exists a non-zero projection $s_0\in  L(\Lambda)'\cap M$ such that  $L(\Lambda)s_0\nprec L(\Gamma).$
\item One shows that there exists a unitary $u\in\mathcal U(M\bar\otimes M)$ such that $$u\Delta(L(\Gamma))u^*\subset L(\Gamma\times\Gamma).$$
\item Next, one proves that the algebra $C:=\Delta(A)'\cap (M\bar\otimes M)$ satisfies $$Cq\prec_{M\bar\otimes M}A\bar\otimes A \text{  for all  } q\in\mathcal Z(C).$$
\item Using the previous steps together with a generalization of \cite[Theorem 5.2]{Po04}, one essentially obtains that there exist a unitary $v\in\mathcal U(M\bar\otimes M)$, a group homomorphism $\delta:\Gamma\to\Gamma\times\Gamma$ and a character $\omega:\Gamma\to\mathbb C$ such that $vCv^*=A\bar\otimes A$ and $v\Delta(u_g)v^*=\omega(g)u_{\delta(g)}$, for all $g\in\Gamma$ (the precise statement is the Step 3 of the proof \cite[Theorem 10.1]{IPV10}).

\item Using Step 5, one proves that for every sequence $(x_n)_n$ in $M$ for which the Fourier coefficient (w.r.t. the decomposition $M=A\rtimes\Gamma$) converges to $0$ pointwise in $\|\cdot\|_2$, then the Fourrier coefficient of $\Delta(x_n)$ (w.r.t. the decomposition $M\bar\otimes M=(M\bar\otimes A)\rtimes\Gamma)$ also converges to $0$ pointwise in $\|\cdot\|_2.$ This shows $B\prec A$ and Step 1 implies that $\Gamma\curvearrowright X$ is W$^*$-superrigid.
\end{enumerate}

\subsection{\bf Proof of Theorem \ref{1}} We record first the following observation.
\begin{remark}\label{free}
Since $\Sigma$ is almost malnormal in $\Gamma$, using \cite[Lemma 5.3]{Dr15}, the action $\Gamma\curvearrowright X$ is free (see also \cite[Lemma 2.1]{Io06b}).
\end{remark}

{\it Proof of Theorem \ref{1}.}
Assume that $\Lambda\curvearrowright (Y,\nu)$ is an arbitrary free ergodic pmp action such that 
\[M:=L^\infty(X)\rtimes\Gamma=L^\infty(Y)\rtimes\Lambda.
\]
We put $A=L^\infty(X), B=L^\infty(Y)$. Define $\Delta:M\to M\bar\otimes M$ by $\Delta(bv_s)=bv_s\otimes v_s,$  for all $b\in B$ and $s\in\Lambda,$ where we denote by $v_s, s\in\Lambda,$ the canonical unitaries corresponding to the action of $\Lambda$. 

Since the action $\Gamma\curvearrowright X$ is OE superrigid (using \cite[Theorem A]{Dr15} and \cite[Theorem 5.6]{Po05}), Step 1 is completed. To prove Step (2), suppose $L(\Lambda)q\prec L(\Gamma)$ for all $q\in L(\Lambda)'\cap M$. Since $\Sigma$ is amenable, $L(\Lambda)\nprec L(\Sigma)$, so by Corollary \ref{unitary}, there exists a unitary $u\in\mathcal U(M)$ such that $uL(\Lambda)u^*\subset L(\Gamma).$  Based on Step 1, Theorem \ref{criterion} proves that $\Gamma\curvearrowright X$ is W$^*$-superrigid. This completes Step 2. Therefore, we take a non-zero projection $q_0\in L(\Lambda)'\cap M$ such that $L(\Lambda)q_0\nprec L(\Gamma).$ Step (3) is obtained by combining Theorem \ref{ipv} and \cite[Lemma 10.2.5]{IPV10}. 

{\it Proof of Step (4).}
Note that Theorem \ref{clustering} proves this step by considering the abelian subalgebra $D_0:=\Delta(A)(1\otimes q_0)$. For showing this, denote $C_0=D_0'\cap (M\bar\otimes q_0Mq_0)$, $C=\Delta(A)'\cap (M\bar\otimes M)$ and note that $C_0=C(1\otimes q_0)$. Since $L(\Lambda)q_0\nprec L(\Gamma)$, \cite[Lemma 9.2.4]{Io10} implies that $\Delta(M)(1\otimes q_0)\nprec M\bar\otimes L(\Gamma).$ Using \cite[Lemma 10.2]{IPV10}, we see that all the conditions of Theorem \ref{clustering} are satisfied. Therefore, we obtain that $C_0q\prec A\bar\otimes A$, for all $q\in \mathcal Z(C_0)=\mathcal Z(C)(1\otimes q_0)$.
\hfill$\square$

{\it Proof of Step (5).} For proving Step (3) of the proof of \cite[Theorem 10.1]{IPV10}, one only needs to show: 
\begin{itemize}
\item If $H$ is a subgroup of $\Gamma\times\Gamma$ such that $H$ acts non-ergodically on $A\bar\otimes A$, 
then $\Delta(L(\Gamma))\nprec (A\bar\otimes A)\rtimes H.$ 
\end{itemize}
Suppose by contrary that $\Delta(L(\Gamma))\nprec (A\bar\otimes A)\rtimes H.$
It is easy to prove that there exists a finite set $T\subset\Gamma$ such that $H\subset (\cup_{t\in T}t\Sigma)\times\Gamma$ or $H\subset \Gamma\times (\cup_{t\in T}t\Sigma).$ This implies that $\Delta(L(\Gamma))\prec (A\rtimes\Sigma)\bar\otimes M$ or $\Delta(L(\Gamma))\prec M\bar\otimes (A\rtimes\Sigma).$ By applying \cite[Lemma 10.2.5]{IPV10}, we obtain a contradiction. \hfill$\square$

Step (6) works in general once the other steps are proven. This finishes the proof of the theorem.

\hfill$\blacksquare$ \\






\subsection{\bf Proof of Theorem \ref{2}}\label{section2} In this subsection we will prove a more general statement of Theorem \ref{2}.

\begin{assumption}\label{assumption2} Let $\Sigma$ be a subgroup of a countable icc group $\Gamma$. Let $\sigma_0$ be a pmp action of $\Sigma$ on a non-trivial standard probability space $(X_0,\mu_0)$ and denote by $\sigma$ the coinduced action of $\Gamma$ on $X:=X_0^{\Gamma/\Sigma}$. Suppose:

\begin{itemize}
\item $\Gamma$ is a non-amenable icc group which is measure equivalent to a group $\Lambda_0$ for which the group von Neumann algebra $L(\Lambda_0)$ is not prime. 

\item $\Sigma$ is amenable and almost malnormal. 
\end{itemize}
\end{assumption}

\begin{theorem}\label{2'}
Let $\Gamma\car X$ be as in Assumption \ref{assumption2}. Then $\Gamma\car X$ is W$^*$-superrigid.
\end{theorem}

{\it Proof.} The proof of this theorem goes along the same lines as the proof of Theorem \ref{1}. We point out only the differences. The action $\Gamma\curvearrowright X$ is OE superrigid using Theorem \ref{cocycle} and \cite[Theorem 5.6]{Po05}. Step (3) follows by Theorem \ref{delta}. All the other steps follow as in the proof of Theorem \ref{1}, which finishes the proof. \hfill$\blacksquare$\\

\begin{remark}\label{remark.}
A careful handling of Thorem \ref{delta} shows that Assumption \ref{assumption2} can be improved by supposing the weaker assumption that $L(\Lambda_0)$ contains a commuting pair of diffuse subalgebras $P_1$ and $P_2$ such that $P_2$ is non-amenable and $\mathcal N_{L(\Lambda_0)}(P_1\vee P_2)''=L(\Lambda_0)$  (see also Step 1 of the proof of \cite[Theorem 8.2]{IPV10}). 
\end{remark}

\begin{corollary}
Let $\Gamma$ be an icc non-amenable group which is measure equivalent to a group $\Lambda_0$ for which $L(\Lambda_0)$ is not prime. Then the Bernoulli action $\Gamma\curvearrowright (X,\mu)^\Gamma$ is $W^*$-superrigid, where $(X,\mu)$ is a non-trivial standard probability space.
\end{corollary}


\begin{thebibliography}{ABC99}

\bibitem [BHV08]{BHV08} B.Bekka, P. de la Harpe, A.Valette: {\it Kazhdan's Property (T)}, (New Mathematical Monographs, 11),
Cambridge University Press, Cambridge, 2008.

\bibitem [Bo63]{Bo63} A. Borel: {\it Compact Clifford-Klein forms of symmetric spaces}, Topology {\bf 2} (1963), 111-122.

\bibitem [Bo12a]{Bo12a} R. Bouttonet: {\it On solid ergodicity for Gaussian actions}, J. Funct. Anal.,{\bf 263} (2012) 1040-1063.
\bibitem [Bo12]{Bo12} R. Boutonnet: {\it W$^*$-superrigidity of mixing Gaussian actions of rigid groups}, Adv. Math. {\bf 244} (2013)
\bibitem [Bo14]{Bo14} R. Bouttonet: {\it Plusieurs aspects de rigidit\'{e} des alg\`{e}bres de von Neumann}, PhD thesis (2014).


\bibitem[BV12]{BV12} M. Berbec, S. Vaes: {\it W$^*$-superrigidity for group von Neumann algebras or left-right wreath products}, Proceedings of the London Mathematical Society {\bf 108} (2014), 1116-1152. 

\bibitem [CI08]{CI08} I. Chifan, A. Ioana: {\it Ergodic subequivalence relations induced by a Bernoulli action}, Geometric and Functional Analysis, 20(1): 53-67, 2010.

\bibitem [CI17]{CI17} I. Chifan, A. Ioana: {\it Amalgamated free product rigidity for group von Neumann algebras}, Preprint, arxiv:1705.07350.
\bibitem [CIK13]{CIK13}  I. Chifan, A. Ioana, Y. Kida:  {\it W$^*$-superrigidity for arbitrary actions of central quotients of braid
groups}, Math. Ann. {\bf 361} (2015), 925-959.
\bibitem [CJ82]{CJ82} A. Connes, V.F.R. Jones: {\it A II$_1$ factor with two non-conjugate Cartan subalgebras},
Bull. Amer. Math. Soc. {\bf 6} (1982), 211-212.

\bibitem [CK15]{CK15} I. Chifan, Y. Kida: {\it OE and W* superrigidity results for actions by surface braid groups}, Proceedings of the London Mathematical Society, in press.

\bibitem [Co76]{Co76} A. Connes: {\it Classification of Injective Factors Cases $II_1, II_\infty, III_\lambda, \lambda \neq 1$}, Ann. of Math., {\bf 104} (1976), 73-115.
\bibitem [Co94]{Co94} A. Connes: {\it Noncommutative Geometry}, Academic Press, 1994.
\bibitem [CP10]{CP10} I. Chifan, J. Peterson: {\it Some unique group-measure space decomposition results}, Duke Math. J. {\bf 162} (2013), no. 11, 1923-1966. 
\bibitem [CS11]{CS11} I. Chifan, T. Sinclair: {\it On the structural theory of II$_1$ factors of negatively curved groups}, Ann. Sci.
Éc. Norm. Supér. (4) {\bf 46} (2013), 1-33.
\bibitem [CSU11]{CSU11} I. Chifan, T. Sinclair, B. Udrea: {\it On the structural theory of II$_1$ factors of negatively curved groups,
II. Actions by product groups}, Adv. Math. {\bf 245} (2013), 208-236.



\bibitem [DHI16]{DHI16} D. Drimbe, D. Hoff, A. Ioana:{\it Prime II$_1$ factors arising from irreducible lattices in products of rank one simple Lie groups} to apear in J. Reine. Angew. Math.

\bibitem [dHW14]{dHW14} P. de la Harpe, C. Weber: {\it Malnormal subgroups and Frobenius groups: basics and
examples}, with an appendix by Denis Osin, Confluentes Math. {\bf 6} (2014), no. 1, 65–76
\bibitem [Dr15]{Dr15} D. Drimbe: {\it Cocycle and orbit equivalence superrigidity for coinduced actions}, to appear in  Ergodic Theory Dynam. Systems.
\bibitem[Fu99]{Fu99} A. Furman: {\it Orbit equivalence rigidity}, Ann. of Math. (2) {\bf 150} (1999), no. 3, 1083-1108.
\bibitem [FV10]{FV10} P. Fima, S. Vaes: {\it HNN extensions and unique group measure space decomposition of II$_1$ factors}
Trans. Amer. Math. Soc. {\bf 354} (2012) 2601-2617.
\bibitem[Ho15]{Ho15} D. Hoff: {\it Von Neumann algebras of equivalence relations with nontrivial one-cohomology}, J. Funct. Anal. {\bf 270} (2016), no. 4, 1501-1536.


\bibitem[GITD16]{GITD16} D. Gaboriau, A. Ioana, R. Tucker-Drob: {\it Cocycle superrigidity for translation actions of product groups}, summited, arXiv 1603.07616 (2016).

\bibitem[Gr91]{Gr91} M. Gromov: {\it Asymptotic invariants of infinite groups, Geometric group theory}, Vol. 2 (Sussex, 1991), London Math. Soc. Lecture Note Ser., vol. 182, Cambridge Univ. Press, Cambridge, 1993, pp. 1-295.

 \bibitem [HPV10]{HPV10} C. Houdayer, S. Popa, S. Vaes: {\it  A class of groups for which every action is W$^*$-superrigid},
Groups Geom. Dyn. {\bf 7} (2013), 577-590.
\bibitem [Io06a]{Io06a} A. Ioana: {\it Rigidity results for wreath product of II$_1$ factors}, J. Funct. Anal. {\bf 252} (2007), 763-791.
\bibitem [Io06b]{Io06b} A. Ioana: {\it Orbit inequivalent actions for groups containing a copy of $F_2$}, Invent. Math. {\bf 185} (2011), 55-73.
\bibitem [Io10]{Io10} A. Ioana: {\it W$^*$-superrigidity for Bernoulli actions of property (T) groups} J. Amer. Math. Soc. {\bf 24} (2011), 1175-1226.

\bibitem [Io12a]{Io12a} A. Ioana: {\it Classification and rigidity for von Neumann algebras}, European Congress of Mathematics, EMS (2013), 601-625.
\bibitem [IPV10]{IPV10} A. Ioana, S. Popa, S. Vaes: {\it A Class of superrigid group von Neumann algebras}, Ann. of Math.
(2) {\bf 178} (2013), 231-286

 
\bibitem [Ka67]{Ka67}  D. Kazhdan: {\it On the connection of the dual space of a group with the structure of
its closed subgroups}, Funct. Anal. and its Appl. {\bf 1}(1967), 63-65. 

 \bibitem[KV15]{KV15} A. Krogager, S. Vaes: {\it A class of II$_1$ factors with exactly two crossed product decompositions}, preprint arXiv:1512.06677.

\bibitem [Ma82]{Ma82}  G. Margulis: {\it Finitely-additive invariant measures on Euclidian spaces}, Ergodic Theory
Dynam. Systems {\bf 2}(1982), 383-396.

\bibitem [MvN36]{MvN36} F. J. Murray, J. von Neumann, On rings of operators. Ann. of Math. 37 (1936).
116-229.

\bibitem[OP07]{OP07} N. Ozawa, S. Popa: {\it On a class of II$_1$ factors with at most one Cartan subalgebra}, Ann. of Math. (2) {\bf 172} (2010), no. 1, 713-749.
\bibitem [Pe09]{Pe09} J. Peterson: {\it Examples of group actions which are virtually W$^*$-superrigid}, Preprint. arXiv:1002.1745.
\bibitem [Po01]{Po01} S. Popa: {\it On a class of type II$_1$ factors with Betti numbers invariants}, Ann. of Math.
{\bf 163} (2006), 809-899.
\bibitem [Po03]{Po03} S. Popa: {\it Strong rigidity of II$_1$ factors arising from malleable actions of w-rigid groups. I.}, Invent. Math. {\bf 165} (2006), 369-408. 
\bibitem [Po04]{Po04} S. Popa: {\it Strong rigidity of II$_1$ factors arising from malleable actions of w-rigid groups. II.}, Invent. Math. {\bf 165}(2) (2006), 409-451.
\bibitem [Po05]{Po05} S. Popa: {\it Cocycle and orbit equivalence superrigidity for malleable actions of {\it w}-rigid groups}, Invent. Math. {\bf 170} (2007), 243-295.
\bibitem [Po06b]{Po06b} S. Popa: {\it On a class of type II$_1$ factors with Betti numbers invariants}, Ann. of Math. {\bf 163} (2006), 809-889.
\bibitem [Po06]{Po06} S. Popa: {\it On the superrigidity of malleable actions with spectral gap}, J. Amer. Math. Soc. {\bf 21} (2008),  981-1000.


\bibitem [Po07]{Po07} S. Popa: {\it Deformation and rigidity for group actions and von Neumann algebras}, In Proceedings of the ICM (Madrid, 2006), Vol. I, European Mathematical Society Publishining House, 2007, 445-477.
\bibitem [PS09]{PS09} J. Peterson, T. Sinclair: {\it On cocycle superrigidity for Gaussian actions}  Erg. Th. and Dyn. Sys.{\bf  32} (2012), no. 1, 249-272.

\bibitem [PV06]{PV06} S. Popa, S. Vaes: {\it Strong rigidity of generalized Bernoulli actions and computations of their symmetry groups}, Adv. Math. {\bf 217} (2008), 833-872.

\bibitem [PV09]{PV09} S.Popa, S. Vaes: {\it Group measure space decomposition of II$_1$ factors and $W^∗$-superrigidity}, Invent. Math. {\bf 182} (2010), no. 2, 371-417.
\bibitem [PV11]{PV11} S. Popa, S. Vaes: {\it Unique Cartan decomposition for II$_1$ factors arising from arbitrary actions of free groups}, Acta Mathematica {\bf 212} (2014), 141-198.
\bibitem [PV12]{PV12} S. Popa, S. Vaes: {\it Unique Cartan decomposition for II$_1$ factors arising from arbitrary actions of hyperbolic groups}, Journal fur die reine und angewandte Mathematik, {\bf 690} 2014, 433-458.

\bibitem [Ra72]{Ra72} M. S. Raghunathan: {\it Discrete subgroups of Lie groups}, Springer-Verlag, New York, 1972, Ergebnisse der Mathematik und ihrer Grenzgebiete, Band 68.

\bibitem [RS10]{RS10} G. Robertson,Tim Steger: {\it Malnormal subgroups of lattices and
the Pukanszky invariant in group factors}, J. Funct. Anal. {\bf 258}-8 (2010), 2708-2713.
\bibitem [Si55]{Si55} I. M. Singer: {\it Automorphism of finite factors}, Amer. J. Math. {\bf 77} (1955), 117-133.
\bibitem[Va08]{Va08} S. Vaes: {\it Explicit computations of all finite index bimodules for a family of II$_1$ factors}, Ann. Sci. \'{E}c. Norm. Sup\'{e}r. (4) {\bf 41} (2008), no. 5, 743-788.
\bibitem [Va10a]{Va10a} S. Vaes: {\it Rigidity for von Neumann algebras and their invariants}, Proceedings of the ICM (Hyderabad, India, 2010), Vol. III, Hindustan Book Agency (2010), 1624-1650.
 \bibitem [Va10b]{Va10b} S Vaes: {\it One-cohomology and the uniqueness of the group measure space decomposition of a II$_1$ factor}, Math. Ann. {\bf 355} (2013), 661-696.




\end{thebibliography}
\end{document}